\documentclass[10pt,a4paper]{amsart}
\usepackage{amssymb}
\usepackage{mathrsfs}  
\usepackage{graphicx}
\usepackage{hyperref}  
\hypersetup{colorlinks=true, linkcolor=blue, anchorcolor=blue, 
citecolor=red, filecolor=blue, menucolor=blue, pagecolor=blue, 
urlcolor=blue}
\numberwithin{equation}{section}
\newtheorem{theorem}{Theorem}[section]

\newtheorem{lemma}[theorem]{Lemma}
\theoremstyle{remark}
\newtheorem{remark}{Remark}[section]

\theoremstyle{definition}
\newtheorem{definition}[theorem]{Definition}
\newtheorem{example}[theorem]{Example}

\newcommand{\xxi}{\langle\xi\rangle}
\newcommand{\xx}{\langle x\rangle}
\newcommand{\DD}{\langle D\rangle}
\newcommand{\DDx}{\langle D_{x}\rangle}
\newcommand{\zz}{\langle z\rangle}
\newcommand{\xxR}{\langle x\rangle_{R} }
\newcommand{\xxS}{\langle x\rangle_{S} }
\newcommand{\DDM}{\langle D\rangle_{M} }
\newcommand{\one}[1]{\mathbf{1}_{#1}}

\begin{document}
\title%
[Non flat waveguides]%
{Evolution equations on non flat waveguides}
\begin{abstract}
  We investigate the dispersive properties of evolution equations
  on waveguides with a non flat shape. More precisely we
  consider an operator
  \begin{equation*}
    H=-\Delta_{x}-\Delta_{y}+V(x,y)
  \end{equation*}
  with Dirichled boundary condition on an unbounded domain $\Omega$,
  and we introduce the notion of a \emph{repulsive waveguide} along
  the direction of the first group of variables $x$. 
  If $\Omega$ is a repulsive waveguide,
  we prove a sharp estimate for the Helmholtz equation $Hu-\lambda u=f$.
  As consequences we prove smoothing estimates for
  the Schr\"odinger and wave equations associated to $H$, and
  Strichartz estimates for the Schr\"odinger equation. Additionally,
  we deduce that the operator $H$ does not admit eigenvalues.
\end{abstract}
\date{\today}    
\author{Piero D'Ancona}
\address{Piero D'Ancona: 
SAPIENZA - Universit\`a di Roma,
Dipartimento di Matematica, 
Piazzale A.~Moro 2, I-00185 Roma, Italy}
\email{dancona@mat.uniroma1.it}

\author{Reinhard Racke}
\address{Reinhard Racke:
Fachbereich Mathematik und Statistik,
Universit\"at Konstanz,
Fach D 187, 78457 Konstanz}
\email{reinhard.racke@uni-konstanz.de}

\maketitle

\section{Introduction}\label{sec:intro} 

\begin{figure}[htbp]
\begin{minipage}[bt]{.4\linewidth}
  \centering
    \includegraphics[height=1in]{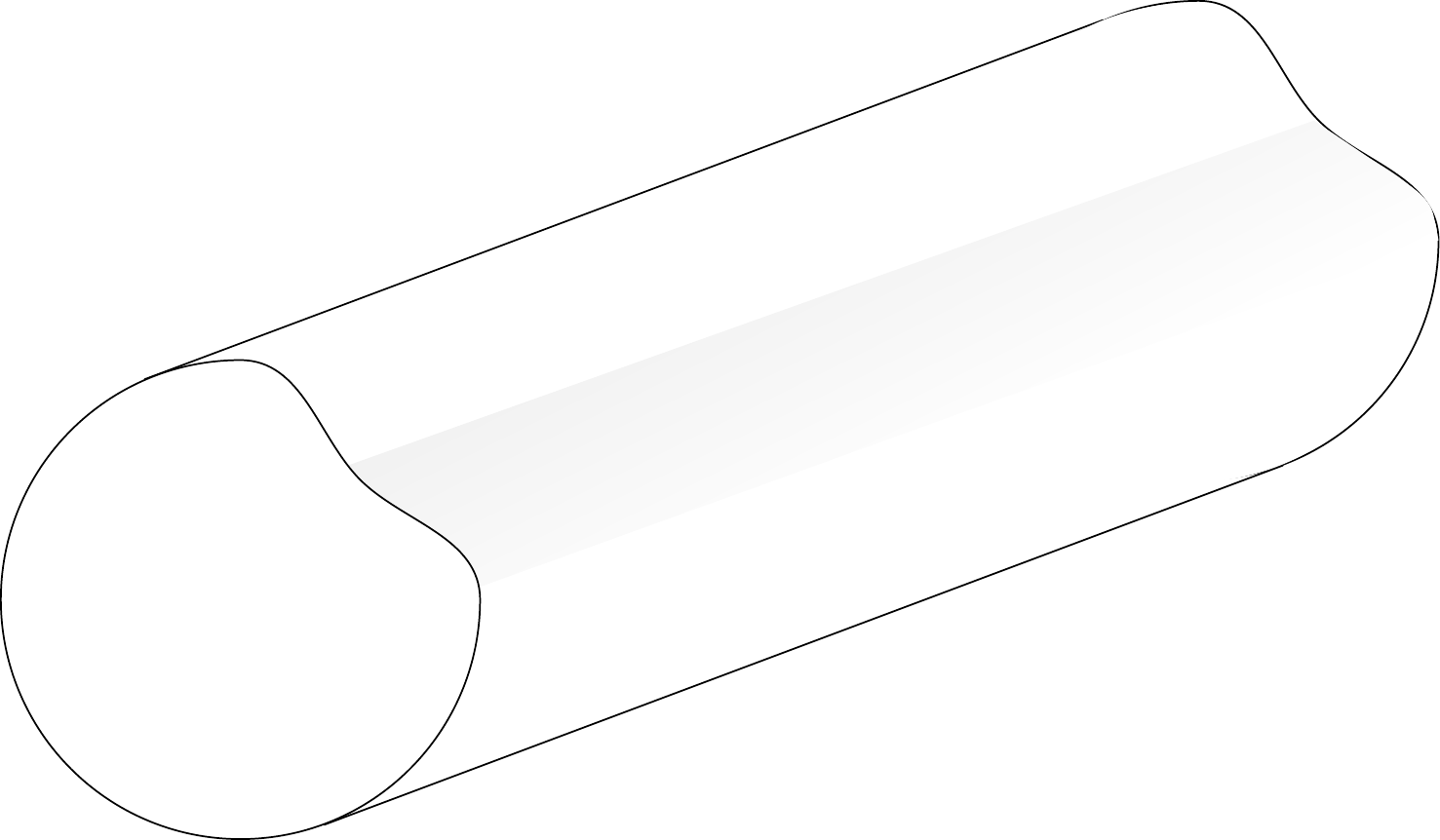} 
\end{minipage}
\hspace{.05\linewidth}
\begin{minipage}[bt]{.4\linewidth}
  \centering
    \includegraphics[height=.6in]{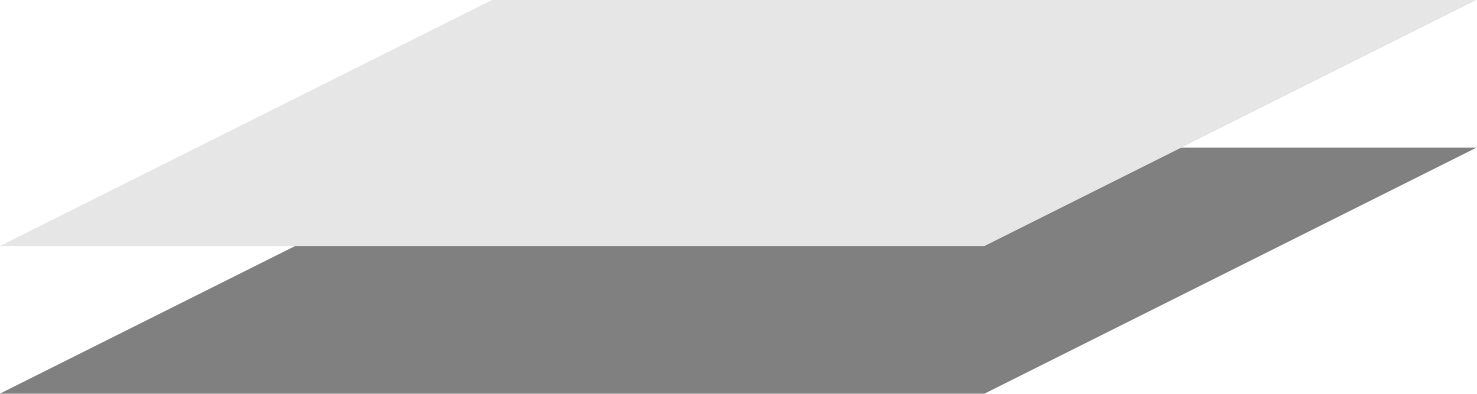}
\end{minipage}
\label{fig:flat}
\caption{(a) $n=1$, $m=2$; (b) $n=2$, $m=1$}
\end{figure}

A \emph{flat waveguide} is a domain $\Omega$ in $\mathbb{R}^{n+m}$ 
which can be written as a product
of a bounded open subset $\omega$ with $\mathbb{R}^{n}$:
\begin{equation*}
  \omega \subseteq \mathbb{R}^{m},\qquad
  \Omega=\mathbb{R}^{n}\times \omega\subseteq \mathbb{R}^{n}_{x}
  \times \mathbb{R}^{m}_{y},\qquad
  n,m\ge1.
\end{equation*}
Throughout the paper we shall denote with $x$ the group of the first
$n$ variables and with $y$ the last $m$ variables in $\mathbb{R}^{n+m}$.
Waveguides appear in many concrete applications,
since they can be used to model various 
interesting physical structures
such as \emph{wires} and \emph{plates} (see Figure \ref{fig:flat}).
The Laplace operator on $\Omega$ with Dirichlet or Neumann boundary
conditions has a natural splitting
\begin{equation*}
  \Delta_{x,y}=\Delta_{x}+\Delta_{y}
\end{equation*}
where $\Delta_{x}$ is the free Laplacian on $\mathbb{R}^{n}$ and
$\Delta_{y}$ is the Dirichlet resp.~Neumann Laplacian on $\Omega$
(we shall also write
\begin{equation*}
  \nabla=(\nabla_{x},\nabla_{y})
\end{equation*}
with obvious meaning).
Thus the operator has a simple spectral structure:
indeed, if we choose an orthonormal set of eigenfunctions
$\{\phi_{j}(y)\}_{j\ge1}$ for $-\Delta_{y}$ on $\omega$ and denote
by $\lambda_{j}^{2}$ the corresponding eigenvalues, the operator
$-\Delta_{x,y}$ is equivalent to the sequence of operators
on $\mathbb{R}^{n}$
\begin{equation*}
  -\Delta_{x}+\lambda_{j}^{2}.
\end{equation*}
As a consequence, the study of linear and
nonlinear evolution equations on
flat waveguides is quite similar to the standard
case of free equations on $\mathbb{R}^{n}$. The theory was
initiated in \cite{LeskyRacke03-a}
and developed in \cite{MetcalfeSoggeStewart05-a} and
\cite{LeskyRacke08-a}. 

Despite the simplicity of the theory, it is clear
that the flatness assumption on the domain is not always realistic.
Thus a natural question is whether a similar theory can be developed for
more general, non flat waveguides. Here we begin to address this
question, by investigating the smoothing and dispersive properties of 
wave and Schr\"odinger equations in more general situations. Such
properties, which are usually expressed as global in time estimates on
solutions of the linear equations, are the key ingredients for 
the nonlinear theory. To the best of our knowledge, the results
in the present paper are the first ones concerning dispersive
phenomena on non flat waveguides.

We start with a quick overview of the dispersive properties for
the linear Schr\"odinger and wave-Klein-Gordon equations in the flat 
case.

\begin{example}\label{exa:schroflat}
  Consider the Schr\"odinger equation
  \begin{equation}\label{eq:schreqhom}
    iu_{t}-\Delta u=0,\qquad
    u(0,x,y)=f(x,y)
  \end{equation}
  with Dirichlet boundary conditions on 
  $\Omega=\mathbb{R}^{n}\times \omega$,
  with $\omega$ a bounded open set in $\mathbb{R}^{m}$. Let
  $\phi_{j}$, $\lambda_{j}^{2}$ be as above, then by expanding
  \begin{equation*}
    u=\sum_{j\ge1}u_{j}(t,x)\phi_{j}(y),
    \qquad
    f=\sum_{j\ge1}f_{j}(x)\phi_{j}(y)
  \end{equation*}
  we can rewrite equation \eqref{eq:schreqhom} as the equivalent family
  of independent equations
  \begin{equation}\label{eq:schreqj}
    i \partial_{t} u_{j}-\Delta_{x} u_{j}+\lambda_{j}^{2} u_{j}
       =0,\qquad
    u_{j}(0,x)=f_{j}(x).
  \end{equation}
  The term $\lambda_{j}^{2}u$ can be absorbed in $iu_{t}$
  via the gauge transformation
  $u_{j}\to e^{i \lambda_{j}^{2}t}u_{j}$, leaving the $L^{p}$
  norm of the solution unchanged. 
  Thus from the explicit representation of the solution
  we have the \emph{dispersive estimates}
  \begin{equation}\label{eq:dispj}
    \|u_{j}(t)\|_{L^{\infty}(\mathbb{R}^{n})}\le |t|^{-n/2}
      \|f_{j}\|_{L^{1}(\mathbb{R}^{n})}
  \end{equation}
  and summing over $j$ we obtain
  \begin{equation}\label{eq:dispsum}
    \|u(t)\|_{L^{\infty}(\Omega)}\le |t|^{-n/2}
    \sum_{j\ge1}\|\phi_{j}\|_{L^{\infty}(\omega)}
    \|f_{j}\|_{L^{1}(\mathbb{R}^{n})}\equiv|t|^{-n/2}\|f\|_{Z}.
  \end{equation}
  A more explicit expression of the norm $\|f\|_{Z}$ requires some
  information on the growth of the maximum norm of eigenfunctions.
  Typically one has
  \begin{equation*}
    \|\phi_{j}\|_{L^{\infty}(\omega)}\lesssim \lambda_{j}^{\sigma}
  \end{equation*}
  for some $\sigma>0$, and this leads to a dispersive estimate 
  of the form
  \begin{equation}\label{eq:dispflat}
    \|u(t)\|_{L^{\infty}(\Omega)}\lesssim
    |t|^{-n/2}\|(1-\Delta_{y})^{\sigma/2+\epsilon} f\|
    _{L^{1}_{x}L^{2}_{y}(\omega)}
  \end{equation}
  
  The pointwise estimate \eqref{eq:dispflat} is quite strong
  and we shall not be able to prove an analogous in the non flat case.
  However Schr\"odinger equations satisfy weaker but more
  general estimates called \emph{Strichartz estimates}, which
  can be extended to our situation. Consider for maximum generality
  the nonhomogeneous equation
  \begin{equation}\label{eq:schreq}
    iu_{t}-\Delta u=F(t,x,y),\qquad
    u(0,x,y)=f(x,y)
  \end{equation}
  with Dirichlet boundary conditions on $\Omega$ as above.
  Here we assume for simplicity $n\ge3$. Expanding again
  \begin{equation*}
    F=\sum_{j\ge1}F_{j}(t,x)\phi_{j}(y)
  \end{equation*}
  we are led to the equations
  \begin{equation}\label{eq:schreqj}
    i \partial_{t} u_{j}-\Delta_{x} u_{j}+\lambda_{j}^{2} u_{j}
       =F_{j}(t,x),\qquad
    u_{j}(0,x)=f_{j}(x).
  \end{equation}
  The endpoint Strichartz estimate (see \cite{GinibreVelo85-d},
  \cite{KeelTao98-a}) for $u_j$ states that
  \begin{equation}\label{eq:strichj}
    \|u_{j}\|_{L^{2}_{t}L^{\frac{2n}{n-2}}_{x}}\lesssim
    \|f_{j}\|_{L^{2}_{x}}+
    \|F_{j}\|_{L^{2}_{t}L^{\frac{2n}{n+2}}_{x}}
  \end{equation}
  with constants independent of $j$. Squaring and summing over $j$
  we obtain the endpoint Strichartz estimate  for
  flat waveguides:
  \begin{equation}\label{eq:strichsum}
    \|u\|_{L^{2}_{t}L^{2}_{y}L^{\frac{2n}{n-2}}_{x}}
    \lesssim
     \|f\|_{L^{2}_{x,y}(\Omega)}+
     \|F\|_{L^{2}_{tL^{2}_{y}}L^{\frac{2n}{n+2}}_{x}}.
  \end{equation}
  We write the estimate in operator form as follows,
  where $\Delta=\Delta_{x,y}$ with Dirichlet b.c.~on $\Omega$,
  $n\ge3$:
  \begin{equation}\label{eq:strichflat}
    \|e^{it \Delta}f\|_{L^{2}_{t}L^{2}_{y}L^{\frac{2n}{n-2}}_{x}}
    \lesssim \|f\|_{L^{2}(\Omega)},\qquad
    \left\|
      \int_{0}^{t}e^{i(t-s)\Delta}F(s)ds
    \right\|_{L^{2}_{t}L^{2}_{y}L^{\frac{2n}{n-2}}_{x}}
    \lesssim
    \|F\|_{L^{2}_{t}L^{2}_{y}L^{\frac{2n}{n+2}}_{x}}.
  \end{equation}
  Similar estimates hold when $n=1,2$.
  
  An even weaker and more general form of estimates are the
  \emph{smoothing estimates}, which go back at least to
  \cite{Kato65-a}, see also \cite{Ben-ArtziKlainerman92-a}.
  For equations \eqref{eq:schreqj} they take the form
  \begin{equation}\label{eq:smoothj}
    \|\xx^{-1/2-\epsilon}|D_{x}|^{1/2}u_{j}\|_{L^{2}_{t}L^{2}_{x}}
    \lesssim\|f_{j}\|_{L^{2}(\mathbb{R}^{n})}
    +\|\xx^{1/2+\epsilon}|D_{x}|^{-1/2}F_{j}\|_{L^{2}_{t}L^{2}_{x}}
  \end{equation}
  where we are using the notations
  \begin{equation*}
    |D_{x}|=(-\Delta_{x})^{1/2},\qquad
    \xx=(1+|x|^{2})^{1/2}.
  \end{equation*}
  Squaring and summing over $j$ we obtain
  \begin{equation}\label{eq:smoosum}
    \|\xx^{-1/2-\epsilon}|D_{x}|^{1/2}u\|_{L^{2}_{t}L^{2}(\Omega)}
    \lesssim \|f\|_{L^{2}(\Omega)}
    +\|\xx^{1/2+\epsilon}|D_{x}|^{-1/2}F\|_{L^{2}_{t}L^{2}(\Omega)}.
  \end{equation}
\end{example}

\begin{example}\label{exa:waveflat}
Consider the wave-Klein-Gordon equation for $u=u(t,x,y)$
\begin{equation}\label{eq:WE}
  u_{tt}-\Delta_{x,y}u+m^{2}u=0, \qquad m\ge0,\qquad(x,y)\in \Omega
  =\mathbb{R}^{n}\times \omega
\end{equation}
with Dirichlet boundary conditions. Proceeding as above
we obtain the family of problems on $\mathbb{R}^{n}$
\begin{equation}\label{eq:WEj}
  \partial^{2}_{t} u_j-\Delta_{x}u_j+
    (\lambda_{j}^{2}+m^{2})u_j=0.
\end{equation}
Notice that in the case of Dirichlet b.c., 
even if we start from a wave equation for $u$ (i.e. $m=0$),
the equations for $u_j$ will always be of Klein-Gordon type
since $\lambda_{j}^{2}>0$ for all $j$. Now, sharp dispersive estimates
are known for the free equations \eqref{eq:WEj}, and summing
over $j$ we shall obtain dispersive estimates for the original equation
\eqref{eq:WE}. Indeed, using the notations
$\DD=(1-\Delta)^{2}$, $\DDM=(M^{2}-\Delta)^{1/2}$, we can represent
the solution of $\square v+M^{2}v=0$ on $\mathbb{R}^{n}_{x}$ as
\begin{equation*}
  v(t,x)=\cos(t\DDM)v(0)+\frac{\sin(t\DDM)}{\DDM}v_{t}(0),
\end{equation*}
thus we see that the solution can be expressed via the
operator $e^{it\DDM}$. To prove a dispersive estimate for it,
we may use the following estimate in terms of Besov spaces
\begin{equation*}
  \|e^{it\DD}f\|_{L^{\infty}_{x}}
    \le \frac{C}{|t|^{n/2}}\|f\|_{B^{\frac n2+1}_{1,1}}
\end{equation*}
(see e.g.~the Appendix of \cite{DanconaFanelli08-a}), 
and by the scaling $v(t,x)\to v(Mt,Mx)$ we obtain
\begin{equation*}
  \|e^{it\DDM}f\|_{L^{\infty}_{x}}\le C\frac{M^{\frac n2}}{|t|^{n/2}}
  \|f(M \cdot)\|_{B^{\frac n2+1}_{1,1}} 
\end{equation*}
The Besov norm in the estimate is not homogeneous, however at least
for $M\ge c_{0}>0$ we get
\begin{equation}\label{eq:dispM}
  |e^{it\DDM}f\|_{L^{\infty}_{x}}
  \le C(c_{0})\frac{M^{n+1}}{|t|^{n/2}}
  \|f\|_{B^{\frac n2+1}_{1,1}}.
\end{equation}
We can now apply this estimate to equation \eqref{eq:WE} 
i.e.~to the sequence of problems \eqref{eq:WEj}.
The relevant operator for \eqref{eq:WE} is
\begin{equation*}
  e^{it(m^{2}-\Delta_{x,y})^{1/2}}f=
  \sum_{j\ge1}e^{it\DD_{M_{j}}}f_{j}(x)\phi_{j}(y),\qquad
  M_{j}^{2}=m^{2}+\lambda_{j}^{2}
\end{equation*}
where of course $f(x,y)=\sum f_{j}(x)\phi_{j}(y)$.
We obtain
\begin{equation*}
  \|e^{it(m^{2}-\Delta_{x,y})^{1/2}}f\|_{L^{\infty}_{x,y}}\le 
  C|t|^{-n/2}\sum_{j\ge1}(m^{2}+\lambda_{j}^{2})^{\frac{n+1}{2}}
  \|f_{j}(x)\|_{B^{\frac n2+1}_{1,1}}\|\phi_{j}\|_{L^{\infty}}.
\end{equation*}
The last sum defines a norm of the initial data $f$ which can
be estimated by the $W^{N,1}$ norm of $f$ for $N$ large enough.
See \cite{LeskyRacke03-a}, \cite{MetcalfeSoggeStewart05-a}
for more details and the applications to
nonlinear wave equations.
Following the same lines, one can prove Strichartz estimates
for the Wave-Klein-Gordon equation on $\Omega$.

Finally, smoothing estimates for the operators
$e^{it\DDM}$ connected to the equation on $\Omega$
\begin{equation*}
  u_{tt}-\Delta_{x,y}u+M^{2}u=0
\end{equation*}
take the form
\begin{equation}\label{eq:smooKGintr}
  \|\xx^{-1/2-\epsilon}e^{it\DDM}f\|_{L^{2}_{t}L^{2}(\Omega)}\lesssim
  \|f\|_{L^{2}(\Omega)}.
\end{equation}

\end{example}

The above approach, based on splitting and diagonalizing 
part of the operator, requires the domain to be of product
type and breaks down for more general domains. 
Even the spectral problem is difficult, as the following
considerations suggest.

\begin{remark}\label{rem:spect}
  For flat waveguides we have
  a purely continuous spectrum, also for {\em certain} locally
  perturbed waveguides, in particular
  for any local perturbation $\Omega$ of
  $(0,1) \times \mathbb{R}^{n-1}$, 
  for which $\nu(x)\cdot x'\leq 0$ holds
  for any $x=(x_1,x')$ on the boundary $\partial\Omega$, see
  construct local perturbations where the Dirichlet Laplacian has
  eigenvalues below its essential spectrum. But there may also exist
  eigenvalues embedded into the essential spectrum; see
  e.g.~\cite{Witsch90-a}, where the following example is constructed.
  Let $D\subset \mathbb{R}^2$ 
  be bounded, star-shaped with respect to
  the origin and invariant under the orthogonal group. Let $\rho\in
  C^0(\mathbb{R}^k)$ 
  be positive, $\rho(x)=1$ for large $|x|$, $\max\,\rho >
  1$. Then the perturbed wave guide
  $$
  \Omega := \cap_{x\in \mathbb{R}^k} \left(\{x\} \times \rho(x)D\right)
  $$
  has an unbounded sequence of multiple eigenvalues embedded into the
  continuous spectrum. Notice that the presence of embedded eigenvalues
  and hence of stationary solutions is in contrast with the decay
  of the solution. Thus we see that suitable conditions of
  \emph{repuslivity} on the shape of the domain are essential
  in order to exclude eigenvalues and ensure dispersion;
  conversely, in presence of bumps in
  the wrong direction, even small, we expect
  in general concentration of energy
  and disruption of dispersion. 
\end{remark}

\begin{figure}[h]
\begin{minipage}{.4\linewidth}
  \centering
    \includegraphics[height=.9in]{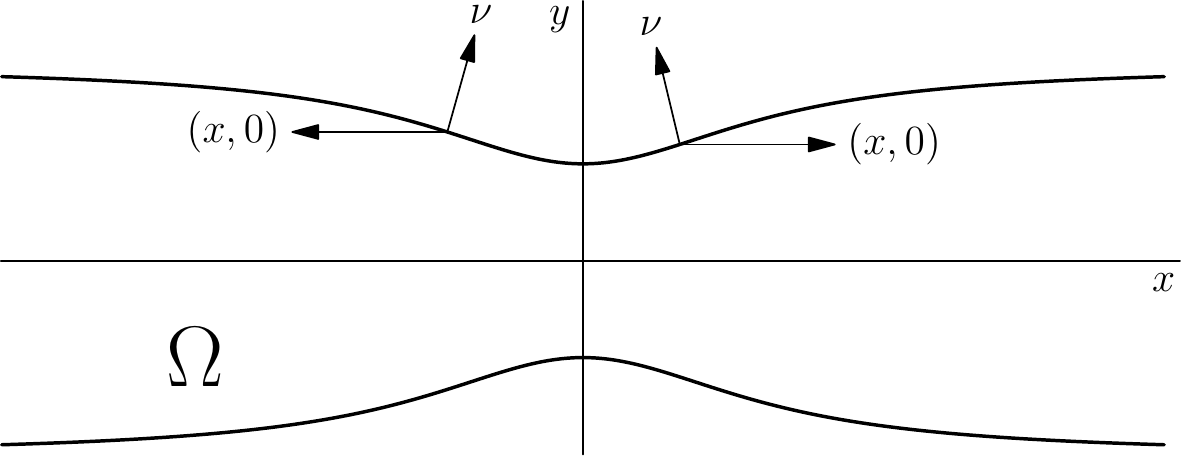}
\end{minipage}
\hspace{.1\linewidth}
\begin{minipage}{.4\linewidth}
  \centering
    \includegraphics[height=1.2in]{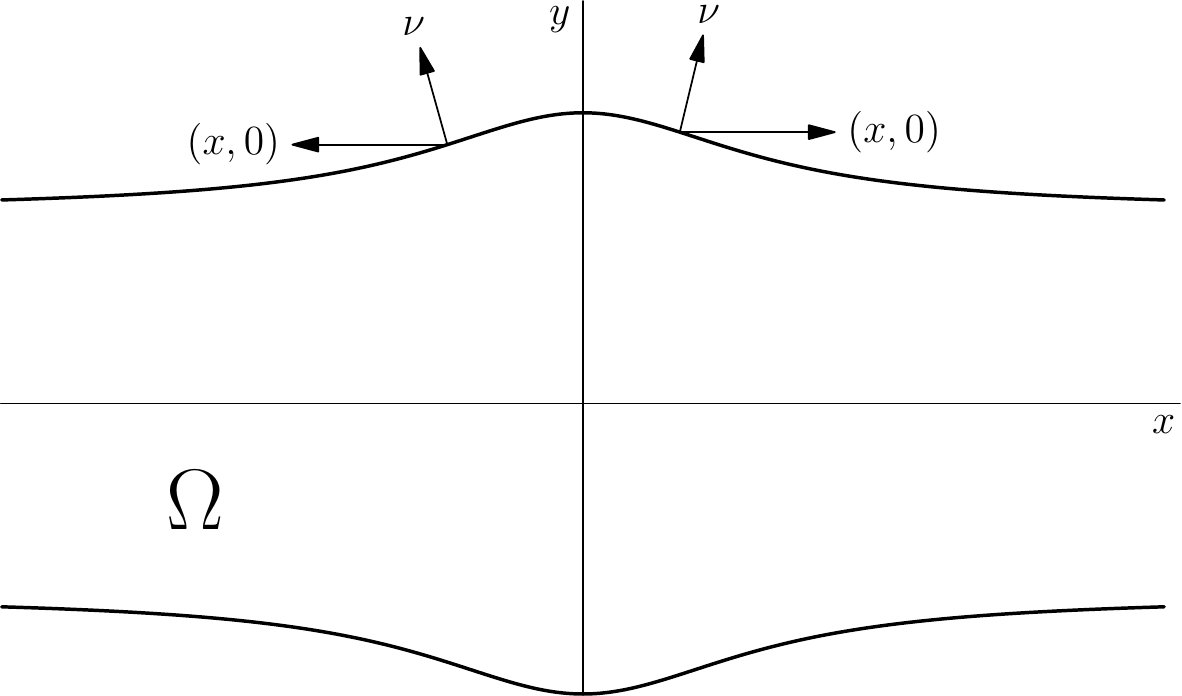}
\end{minipage}
\caption{A repulsive (left) and nonrepulsive (right) domain w.r.to $x$}
\label{fig:rep}
\end{figure}

In order to ensure dispersion, it is reasonable to assume
that the sections of $\Omega$ at fixed $y$
\begin{equation*}
  \{x\in \mathbb{R}^{n}\colon (x,y)\in \Omega\}
\end{equation*}
be nontrapping exterior domains. Actually, in order to prove
smoothing we shall need the following stronger condition
(see Figure \ref{fig:rep}): 

\begin{definition}\label{def:rep}
  Let $\Omega$ be an open subset of 
  $\mathbb{R}^{n}_{x}\times \mathbb{R}^{m}_{y}$
  with Lipschitz boundary, $n,m\ge1$.
  We say that $\Omega$ is \emph{repulsive with respect to the $x$
  variables} if, denoting by $\nu$ the exterior normal to 
  $\partial \Omega$, we have at all points of the boundary
  \begin{equation}\label{eq:repulsive}
    \nu \cdot(x,0)\le0.
  \end{equation}
\end{definition}

We can now state our results. We shall always consider a
waveguide $\Omega$ satisfying condition \eqref{eq:repulsive}, 
with $n\ge3$ and $m\ge1$, and a selfadjoint
Schr\"odinger operator
\begin{equation*}
  H= -\Delta u +V(x,y)
\end{equation*}
with Dirichlet b.c.,
with a locally bounded potential $V(x,y)$ satisfying the assumptions
\begin{equation}\label{eq:assVintro}
  V\ge0,\qquad -x \cdot \nabla_{x}(|x|V)\ge0.
\end{equation}
The conditions on the potential can be substantially relaxed, for
instance by admitting a negative part, small in a suitable sense. 
We did not strive for maximum generality.

\subsection*{Resolvent estimate}\label{sub:resolvent_estimate}
Our approach is based on the Kato smoothing theory 
(see \cite{Kato65-a}, see also
\cite{RodnianskiSchlag04-a}). The crucial tool, which can
be considered the fundamental result of the paper, is a
uniform resolvent estimate for the operator $H$.
To this end we adapt the method of Morawetz multipliers in the 
version of \cite{BarceloRuizVega06-a}. Using the non isotropic
Morrey-Campanato norms
\begin{equation*}%
  \|f\|_{X}=\sup_{R>0}R^{-1/2}\|f\|_{L^{2}(|x|\le R)},\qquad
  \|f\|_{X_{1}}=\sup_{R>0}R^{-3/2}\|f\|_{L^{2}(|x|\le R)},
\end{equation*}
\begin{equation*}
  \|f\|_{X^{*}}=\sum_{j\in \mathbb{Z}}2^{j/2}
      \|f\|_{L^{2}(2^{j-1}\le |x|\le 2^{j})}
\end{equation*}
(which are asymmetric in $x$ and $y$),
our estimate for the resolvent operator $R(z)=(H-z)^{-1}$
can be stated as follows
\begin{equation*}
  \|\nabla_{x}R(z)f\|_{X}^{2}+\|R(z)f\|_{X_{1}}^{2}
  +|z|\|R(z)f\|_{X}^{2}
  \le
  5000n^{2}\|f\|_{X^{*}}^{2}
\end{equation*}
for all $z\not\in \mathbb{R}$ (see Theorem \ref{the:resest}).

\subsection*{Smoothing estimates}\label{sub:smoothing_estimates}  
Using the previous resolvent estimate, an application
of Kato's theory of smooth operators 
allows us to prove the following
smoothing estimates for the Schr\"odinger flow $e^{itH}$
\begin{equation}\label{eq:smoHV1intro}
  \|\xx^{-1/2-\epsilon} |D_{x}|^{1/2} e^{itH}f\|_{L^{2}_{t}L^{2}(\Omega)}+
  \|\xx^{-1-\epsilon}e^{itH}f\|_{L^{2}_{t}L^{2}(\Omega)}\lesssim
  \|f\|_{L^{2}(\Omega)},
\end{equation}
while the nonhomogeneous form of the estimates is
\begin{equation}\label{eq:smoointro}
\begin{split}
  \left\|\xx^{-1/2-\epsilon}
    \int_{0}^{t}  \nabla_{x} e^{i(t-s)H}F(s)ds
  \right\|_{L^{2}_{t}L^{2}(\Omega)}+&
  \\
  +\left\|\xx^{-1-\epsilon}
    \int_{0}^{t} e^{i(t-s)H}F(s)ds
  \right\|_{L^{2}_{t}L^{2}(\Omega)}
  &\lesssim \|\xx^{1+\epsilon} F\|_{L^{2}_{t}L^{2}(\Omega)}
\end{split}
\end{equation}
(see Theorems \ref{the:smoosch}, \ref{the:smoosch2},
\ref{the:smoVH3bis}). 
On the other hand, for the wave-Klein-Gordon equation we prove 
the estimate ($\mu\ge0$)
\begin{equation}\label{eq:smowaveintro}
  \|\xx^{-1/2-\epsilon}e^{it \sqrt{H+\mu^{2}}}f\|_{L^{2}_{t}L^{2}(\Omega)}
  \lesssim
  \|f\|_{L^{2}(\Omega)}
\end{equation}
and, for the inhomogeneous operator,
\begin{equation}\label{eq:smowaveinhintro}
  \left\|\int_{0}^{t}\xx^{-1/2-\epsilon}e^{i(t-s) \sqrt{H+\mu^{2}}}F(s)ds
  \right\|_{L^{2}_{t}L^{2}(\Omega)}
  \lesssim \|\xx^{1/2+\epsilon} F\|_{L^{2}_{t}L^{2}(\Omega)}
\end{equation}
(see Theorem \ref{the:smoowave}). Notice that our results are
comparable with the flat case outlined in
Examples \ref{exa:schroflat} and \ref{exa:waveflat}.

\subsection*{Strichartz estimates}\label{sub:strichart_estimates}  
A typical application of the smoothing estimates is to deduce
Strichartz estimates. We were only able to prove
Strichartz estimates for the Schr\"odinger flow $e^{itH}$,
under the additional assumption that the waveguide 
$\Omega$ coincides with
a flat waveguide outside some bounded region. In this case,
we can recover the full set of Strichartz estimates, however
with a loss of 1/2 derivatives: indeed,
we can prove for all $n\ge3$ and $m\ge1$ the endpoint estimate
\begin{equation}\label{eq:strichnintro}
  \|e^{itH}f\|_{L^{2}_{t}L^{2}_{y}L^{\frac{2n}{n-2}}_{x}}
  \lesssim
  (1+\|\xx^{1+\epsilon} V\|_{L^{2}_{y}L^{n}_{x}})
  \Bigl(
    \|f\|_{L^{2}(\Omega)}+\||D_{x}|^{1/2}f\|_{L^{2}(\Omega)}
  \Bigr)
\end{equation}
(see Theorem \ref{the:strichschro}).

\subsection*{Absence of eigenvalues}\label{sub:eigenvalues}  

As an immediate corollary to the smoothing estimates, we
deduce that, under the conditions on the domain and on the
Schr\"odinger operator $H$ given above (i.e., $n\ge3$,
$m\ge1$, $\Omega$ repulsive w.r.to $x$ and $V$ as in
\eqref{eq:assVintro}),
there are no eigenvalues of $H$, since the presence
of bound states would contradict the $L^{2}$ integrability
in time of the solution. This
generalizes the known results for the special cases in
\cite{Faulhaber82-a} and \cite{MorgenrotherWerner87-a}
described in Remark \ref{rem:spect}.

\bigskip

The natural domain of application of our estimates are
problems of local and global existence for nonlinear evolution
equations. We prefer not to pursue this line of research here;
the applications to nonlinear Schr\"odinger and wave equations 
on non flat waveguides will be the object of future works.


\section{A resolvent estimate}\label{sec:resolvent}  

This section is devoted to a study of the resolvent equation
$u=R(\lambda+i \epsilon)f$ or equivalently
\begin{equation}\label{eq:reseq}
  -\Delta u -(\lambda+i \epsilon)u+V(x,y)u=f.
\end{equation}
We shall follow the classical Morawetz multiplier method
\cite{Morawetz68-a}, in the framework of Morrey-Campanato
spaces as introduced in \cite{PerthameVega99-a},
see also \cite{BarceloRuizVegaVilela-a} and \cite{DanconaFanelli08-a}.
Here additional difficulties are the presence of a boundary,
and the necessity to handle the variables $x$ and $y$ in a
different way.
Moreover, our estimate \eqref{eq:fundresest} is stronger
than the results
in \cite{BarceloRuizVegaVilela-a} in that it provides a uniform
control of the operator $\xx^{-1/2-}|z|^{1/2}R(z)\xx^{-1/2-}$
(corresponding to the last term at the l.h.s. of \eqref{eq:fundresest});
this will allow us to prove a sharp smoothing estimate for the
wave equation in Theorem \ref{the:smoowave}.

The Morrey-Campanato type norms needed here are the following:
\begin{equation}\label{eq:MC}
  \|f\|_{X}=\sup_{R>0}R^{-1/2}\|f\|_{L^{2}(|x|\le R)},\qquad
  \|f\|_{X_{1}}=\sup_{R>0}R^{-3/2}\|f\|_{L^{2}(|x|\le R)},
\end{equation}
\begin{equation}\label{eq:Xs}
  \|f\|_{X^{*}}=\sum_{j\in \mathbb{Z}}2^{j/2}
      \|f\|_{L^{2}(2^{j-1}\le |x|\le 2^{j})}
\end{equation}
and
\begin{equation}\label{eq:surfMC}
  \|f\|_{X_{2}}=\sup_{R>0}R^{-1}\|f\|_{L^{2}(|x|=R)}.
\end{equation}
Notice that the decomposition involves the variables $x$ only.
The $X^{*}$ norm is actually dual to the $X$ norm,
but we shall not need this fact.
For functions $f\in L^{2}_{loc}(\Omega)$ we extend the definition
of these norms by restriction, meaning that
\begin{equation*}
  \|f\|_{X}=\|Ef\|_{X}, \qquad Ef=f \ \text{on $\Omega$,}
  \qquad Ef=0 \ \text{on $\mathbb{R}^{n}\setminus\Omega$,}
\end{equation*}
We shall use the following elementary inequalities:
\begin{equation}\label{eq:MCin1}
  \|fg\|_{L^{1}(\Omega)}\le \|f\|_{X}\|g\|_{X^{*}},
\end{equation}
\begin{equation}\label{eq:MCin3}
  \|fg\|_{L^{1}(\Omega\cap\{R\le|x|\le 2R\})}
  \le 4R^{2}\|f\|_{X}\|g\|_{X_{1}}
\end{equation}
and
\begin{equation}\label{eq:MCin4}
  \|fgh\|_{L^{1}(\Omega)}
  \le 2\|f\|_{X_{1}}\|g\|_{X^{*}}\||x|h\|_{L^{\infty}}
\end{equation}
which implies in particular
\begin{equation}\label{eq:MCin2}
  \|fg\|_{L^{1}(\Omega\cap\{|x|\le R\})}
    \le 2R \|f\|_{X_{1}}\|g\|_{X^{*}},
\end{equation}
Moreover it is easy to see that
\begin{equation}\label{eq:comparnorm}
  \|f\|_{X_{1}}\le \|f\|_{X_{2}}.
\end{equation}

It will also be useful in the following to compare the above norms
with standard weighted $L^{2}$ norms, with weights of the form
\begin{equation}\label{eq:xxR}
  \xxR=(R+|x|^{2}/R)^{1/2},\qquad \xx=(1+|x|^{2})^{1/2}.
\end{equation}
We notice that for all real $s>0$, and for $u$ defined on $\Omega$ 
(after extending
$u$ as zero on $\mathbb{R}^{n}\times \mathbb{R}^{m}$ outside $\Omega$
for simplicity of notation)
\begin{equation*}
  \begin{split}
  \int(R+|x|^{2}/R)^{-s}|u|^{2}dxdy\le
  R^{-s}\int_{|x|\le R}|u|^{2}+
  R^{s}\int_{|x|> R}|x|^{-2s}|u|^{2}
  \\
  \le 
  R^{-s}\int_{|x|\le R}|u|^{2}+
  2^{2s}\sum_{j\ge j_{R}}R^{s} 2^{-2js}
     \int_{C_{j}}|u|^{2}
  \end{split}
\end{equation*}
where $j_{R}=[\log_{2}R]$ and 
$C_{j}=\left\{(x,y):2^{j-1}\le|x|<2^{j}\right\}$. The second term is
bounded by
\begin{equation*}
  \left(
  \sup_{\rho>0}\rho^{-s}\int_{|x|<\rho}|u|^{2}
  \right)
  2^{2s}R^{s}\sum_{j\ge j_{R}}2^{-js}
  \le
  \frac{2^{2s}}{1-2^{-s}}
  \left(
  \sup_{\rho>0}\rho^{-s}\int_{|x|<\rho}|u|^{2}
  \right)
\end{equation*}
so that we have the inequality
\begin{equation}\label{eq:MCtoweightgen}
  \int\xxR^{-2s}|u|^{2}dxdy\le
  \frac{2^{4s}}{2^{s}-1}
  \sup_{\rho>0}\frac{1}{\rho^{s}}\int_{|x|<\rho}|u|^{2}.
\end{equation}
In particular we have, for any $R>0$,
\begin{equation}\label{eq:MCtoweight}
  \|\xxR^{-1}u\|_{L^{2}(\Omega)}\le
  4\|u\|_{X},\qquad
  \|\xxR^{-3}u\|_{L^{2}(\Omega)}\le
  10\|u\|_{X_{1}}.
\end{equation}
By a similar proof we obtain for any $R>0$
\begin{equation}\label{eq:weighttoMC}
  \|u\|_{X^{*}}\le 16
  \|\xxR u\|_{L^{2}(\Omega)}.
\end{equation}
Finally, we notice the following inequality, valid for
all $\gamma>0$ and $\epsilon>0$:
\begin{equation}\label{eq:weight1}
  \|\xx^{-\frac{\gamma}{2}-\epsilon}u\|_{L^{2}}\le
  C(\gamma,\epsilon)
  \sup_{R>0}
  \|\xxR^{-\gamma}u\|_{L^{2}}.
\end{equation}
which evidently holds also with $L^{2}(\Omega)$ in place
of $L^{2}$. To prove it is sufficient to write
\begin{equation*}
  \int\xx^{-\gamma-2 \epsilon}|u|^{2}\le
  \int_{|x|\le1}|u|^{2}+\sum_{j\ge0}2^{-j(\gamma+2 \epsilon)}
     \int_{2^{j}\le|x|<2^{j+1}}|u|^{2}
\end{equation*}
\begin{equation*}
  \le (1+2^{\gamma})\sum_{j\ge0}2^{-2j \epsilon}
  \sup_{R>0}\frac{1}{R^{\gamma}}\int_{|x|\le R}|u|^{2}
\end{equation*}
and observe that
\begin{equation*}
  \frac{1}{R^{\gamma}}\one{|x|\le R}\le
  2^{\gamma}
  \xxR^{-2\gamma}.
\end{equation*}

\begin{theorem}\label{the:resest}
  Let $\Omega \subseteq \mathbb{R}^{n}_{x}\times \mathbb{R}^{m}_{y}$,
  $n\ge3$, $m\ge1$, be a domain repulsive with respect to the
  variables $x$, with Lipschitz boundary. 
  Assume the potential $V(x,y)$ satisfies
  \begin{equation}\label{eq:assV}
    V\ge0,\qquad -\partial_{x}(|x|V)\ge0
  \end{equation}
  and let $u(x,y)\in H^{1}_{0}(\Omega)$ be a solution of
  equation \eqref{eq:reseq}. Then the following estimate holds:
  \begin{equation}\label{eq:fundresest}
    \|\nabla_{x}u\|_{X}^{2}+\|u\|_{X_{1}}^{2}
    +(|\lambda|+|\epsilon|)\|u\|_{X}^{2}
    \le
    5000n^{2}\|f\|_{X^{*}}^{2}.
  \end{equation}
\end{theorem}

\begin{proof}
Consider two real valued
functions $\psi(x)$ and $\phi(x)$,
\emph{independent of the variable} $y$, such that
\begin{equation}\label{eq:assphipsi}
  \nabla \psi, \Delta \psi, \nabla \Delta \psi,\phi,\nabla \phi
  \ \text{are bounded for $|x|$ large}.
\end{equation}
and
\begin{equation}\label{eq:asspsi}
  \nu \cdot \nabla \psi \le0 \ \text{at $\partial \Omega$}.
\end{equation}
Notice that for a function
$\psi(x)$ depending only on $x$ in a radial way, we have
\begin{equation*}
  \nu \cdot\nabla \psi=\nu \cdot(x,0)|x|^{-1}\partial_{x}\psi
\end{equation*}
and recalling Definition \ref{def:rep},
we see that \eqref{eq:asspsi} is equivalent to the condition
that the radial derivative of $\psi$ be non negative:
\begin{equation}\label{eq:asspsi2}
  x \cdot \nabla_{x}\psi\ge0.
\end{equation}
Then we can form the Morawetz multiplier
\begin{equation}\label{eq:mult}
  (\Delta \psi-\phi)\overline{u}+2 \nabla \psi \cdot \nabla \overline{u}.
\end{equation}
Multiplying the resolvent
equation \eqref{eq:reseq} by the quantity \eqref{eq:mult}
and taking the real part we obtain the identity
\begin{equation}\label{eq:ident1}
  \begin{split}
    \nabla u(2D^{2}\psi-\phi I)
    &
    \nabla \overline{u}+
    \frac12 \Delta(\phi-\Delta \psi)|u|^{2}
    +\phi \lambda|u|^{2}
    -(\nabla V \cdot \nabla \psi+\phi V)|u|^{2}
    +\nabla \cdot \Re Q_{1}
    =
       \\
    &
    =\nabla \cdot\Re Q+
    \Re f(2 \nabla \psi \cdot \nabla \overline{u}
              +(\Delta \psi-\phi) \overline{u})
      - 2 \epsilon\Im(\nabla \psi \cdot \nabla \overline{u}\ u)
  \end{split}
\end{equation}
where
\begin{equation}\label{eq:Q}
  Q =\Delta \psi \overline{u} \nabla u 
     -\frac12 \nabla \Delta \psi|u|^{2}
     -(V-\lambda)\nabla \psi|u|^{2}
     +\frac12 \nabla \phi|u|^{2}
     -\phi \overline{u}\nabla u
\end{equation}
and
\begin{equation}\label{eq:Q1}
  Q_{1}=\nabla \psi|\nabla u|^{2}
  -2\nabla u(\nabla \psi \cdot 
  \nabla\overline{u})
\end{equation}

Our goal is to integrate \eqref{eq:fundest} on $\Omega$, with a suitable
choice of the weights $\phi$ and $\psi$.
First of all we show how to handle the last term at the right hand side.
Multiplying \eqref{eq:reseq} by $\overline{u}$
and splitting real and imaginary parts we obtain the two identities
\begin{equation}\label{eq:imu}
  \Im \nabla \cdot\left\{\nabla u \overline{u}\right\}
  +\epsilon|u|^{2}=-\Im(f \overline{u})
  \quad \implies \quad
  \pm\Im \nabla \cdot\left\{\nabla u \overline{u}\right\}
  +|\epsilon||u|^{2}=\mp\Im(f \overline{u})
\end{equation}
$\pm$ being the sign of $\epsilon$, and
\begin{equation}\label{eq:reu}
  \Re \nabla \cdot\left\{-\nabla u \overline{u}\right\}
  +|\nabla u|^{2}=(\lambda-V)|u|^{2}+\Re(f \overline{u}).
\end{equation}
From the second one we deduce (with $\lambda^{+}=\max\{\lambda,0\}$)
\begin{equation*}
  |\epsilon||\nabla u|^{2}\le
  |\epsilon| \lambda^{+}|u|^{2}+|\epsilon|\Re(f \overline{u})+
        \Re \nabla \cdot\left\{|\epsilon| \nabla u \overline{u}\right\}
\end{equation*}
by the positivity of $V(x,y)$, and using \eqref{eq:imu}
\begin{equation*}
  =\mp\lambda^{+}\Im(f \overline{u})
    +\nabla \cdot\left\{\pm\Im \lambda^{+}u\nabla \overline{u}
       +\Re |\epsilon| \nabla u \overline{u}
    \right\}
    +|\epsilon|\Re(f \overline{u})
\end{equation*}
and hence
\begin{equation}\label{eq:reu2}
  |\epsilon||\nabla u|^{2}\le
    (\lambda^{+}+|\epsilon|)|f \overline{u}|
    +\nabla \cdot\left\{\pm\Im \lambda^{+}u\nabla \overline{u}
       +\Re |\epsilon| \nabla u \overline{u}
    \right\}.
\end{equation}
Now by Cauchy-Schwarz we can write
\begin{equation*}
  2 |\epsilon u\nabla \overline{u}|\le
      |\epsilon|(\lambda^{+}+|\epsilon|)^{1/2}|u|^{2}+
      |\epsilon|(\lambda^{+}+|\epsilon|)^{-1/2}|\nabla u|^{2}
\end{equation*}
and using \eqref{eq:imu}, \eqref{eq:reu2}
\begin{equation*}
  \le
  2(\lambda^{+}+|\epsilon|)^{1/2}|f \overline{u}|\mp
  \nabla \cdot\left\{\Im\nabla \overline{u}u\right\}
           (\lambda^{+}+|\epsilon|)^{1/2}
  + \nabla \cdot\left\{\pm\Im \lambda^{+}u\nabla \overline{u}
              +\Re |\epsilon| \nabla u \overline{u}\right\}
           (\lambda^{+}+|\epsilon|)^{-1/2}.
\end{equation*}
In conclusion we have the estimate
\begin{equation}\label{eq:auxest}
  2 |\epsilon u\nabla \overline{u}|\le
  2(\sqrt{|\epsilon|} +\sqrt{\lambda^{+}})|f \overline{u}|
  + \nabla \cdot A
\end{equation}
with
\begin{equation}\label{eq:A}
  A=\frac{|\epsilon|\Re \nabla u \overline{u}
          \pm(2 \lambda^{+}+|\epsilon|)\Im \nabla \overline{u}u}
           {(\lambda^{+}+|\epsilon|)^{1/2}},
    \qquad
    \pm=\ \text{sign of $\epsilon$}.
\end{equation}
We insert this in our basic identity \eqref{eq:ident1} obtaining
the inequality
\begin{equation}\label{eq:fundest}
  \begin{split}
    \nabla u(2D^{2}\psi-\phi I) & \nabla \overline{u}+
    \frac12 \Delta(\phi-\Delta \psi)|u|^{2}
    +\phi \lambda|u|^{2}
    -(\nabla V \cdot \nabla \psi+\phi V)|u|^{2}
    +\nabla \cdot\Re Q_{1}
    \le
       \\
    &
    \le
    2|f \nabla \psi \cdot \nabla \overline{u}|
    +|f(\Delta \psi-\phi) \overline{u}|
    +2\|\nabla \psi\|_{L^{\infty}}
    (\sqrt{|\epsilon|}+ \sqrt{\lambda^{+}})|f \overline{u}|
    +\nabla \cdot \Re P
  \end{split}  
\end{equation}
where
\begin{equation}\label{eq:P}
  P=Q+\|\nabla \psi\|_{L^{\infty}}A
\end{equation}
with $A,Q,Q_{1}$ given by \eqref{eq:A}, \eqref{eq:Q} 
\eqref{eq:Q1} respectively.

Next we show how to estimate the integral over $\Omega$ of the right
hand side of \eqref{eq:fundest}. We need an additional estimate, 
obtained by multiplying \eqref{eq:reseq}
by $\chi \overline{u}$ and taking the imaginary part: as in \eqref{eq:imu}
we get
\begin{equation}\label{eq:imu2}
  \pm\Im \nabla \cdot\left\{\chi\nabla u \overline{u}\right\}
  +|\epsilon|\chi|u|^{2}=\mp\Im(\chi f \overline{u})
  \mp\Im(\nabla \chi \cdot \nabla \overline{u}u).
\end{equation}
We choose $\chi$ as a radial function of the variables $x$ only, and
precisely
\begin{equation*}
  \chi=
  \begin{cases}
    1 &\text{if $ |x|<R $,}\\
    0 &\text{if $ |x|>2R $,}\\
    2-|x|/R &\text{if $ R\le |x|\le 2R $.}
  \end{cases}
\end{equation*}
Then integrating \eqref{eq:imu2}
on $\Omega$ and noticing that the boundary 
terms disappear (thanks to the Dirichlet b.c.), 
we arrive at the inequality
\begin{equation*}
  |\epsilon|\int_{\Omega\cap\{|x|\le R\}}|u|^{2}
  \le\int_{\Omega\cap\{|x|\le 2R\}}|f \overline{u}|+
  \frac1R
  \int_{\Omega\cap\{R\le|x|\le 2R\}}|\nabla_{x} u| |u|
\end{equation*}
since $\chi$ depends only on $x$.
We estimate the right hand side using \eqref{eq:MCin2}, \eqref{eq:MCin3},
and dividing by $R$ we obtain
\begin{equation*}
  \frac{|\epsilon|}R
  \int_{\Omega\cap\{|x|\le R\}}|u|^{2}\le
  4\|f\|_{X^{*}}\|u\|_{X_{1}}
  +4\|\nabla_{x} u\|_{X}\|u\|_{X_{1}}
\end{equation*}
and taking the sup in $R$ we conclude
\begin{equation}\label{eq:basic2}
  |\epsilon|\|u\|_{X}^{2}\le
  4\left(\|f\|_{X^{*}}+\|\nabla_{x} u\|_{X}\right) \|u\|_{X_{1}}
\end{equation}
Now consider the quantity
\begin{equation*}
   2(\sqrt{\lambda^{+}}+\sqrt{|\epsilon|} )
  \|f \overline{u}\|_{L^{1}(\Omega)}\le
  2(\sqrt{\lambda^{+}}+\sqrt{|\epsilon|} )
  \|f\|_{X^{*}}\|u\|_{X}
\end{equation*}
where we used again \eqref{eq:MCin1}.
By \eqref{eq:basic2} we have
\begin{equation*}
  \le
  2\sqrt{\lambda^{+}}\|f\|_{X^{*}}\|u\|_{X}+
  4\|f\|_{X^{*}}(\|f\|_{X^{*}}+\|\nabla_{x}u\|_{X})^{1/2}
      \|u\|_{X_{1}}^{1/2}
\end{equation*}
and hence, for all $\delta\in(0,1)$,
\begin{equation}\label{eq:basic3}
  2(\sqrt{\lambda^{+}}+\sqrt{|\epsilon|} )
  \|f \overline{u}\|_{L^{1}(\Omega)}
  \le\delta(\lambda^{+}\|u\|_{X}^{2} 
       +\|\nabla_{x} u\|_{X}^{2}+\|u\|_{X_{1}}^{2})
   +5 \delta^{-1}\|f\|_{X^{*}}^{2}.
\end{equation}
This inequality will be used to estimate the third term in the
r.h.s. of \eqref{eq:fundest}.

We consider now the term 
$\nabla \cdot \Re P
=\Re\nabla \cdot(Q+\|\nabla \psi\|_{L^{\infty}}A)$, which
vanishes after integration. To see this, we define the cylinder
\begin{equation*}
  C_{R}=\left\{(x,y)\colon
    |x|<R,\ y\in \mathbb{R}^{m}
  \right\},
\end{equation*}
we integrate $\nabla \cdot P$ on 
$\Omega\cap C_{R}$ and let $R\to+\infty$. The boundary of
$\Omega\cap C_{R}$ is the union of the two sets
\begin{equation*}
  S_{1}=\partial\Omega\cap C_{R}
  \quad\text{and}\quad
  S_{2}=\partial C_{R}\cap \Omega=
  \left\{(x,y)\in \Omega \colon |x|=R\right\}
\end{equation*}
and orrespondingly, we get two surface integrals. The integral on
$S_{1}$ vanishes thanks to the Dirichlet boundary condition,
thus we are left with the boundary integral
\begin{equation*}
  \int_{S_{2}}\nu \cdot P d \sigma.
\end{equation*}
By the first assumption \eqref{eq:assphipsi} on the weights $\phi,\psi$
we have evidently
\begin{equation}\label{eq:liminf}
  \liminf_{R\to+\infty}\int_{S_{2}}\nu \cdot P d \sigma=0
\end{equation}
since the function $u$ is in $H^{1}(\Omega)$. This proves that
\begin{equation*}
  \int_{\Omega}(\nabla \cdot P) dx dy=0.
\end{equation*}

Concerning the first and the second term at the right hand side
of \eqref{eq:fundest}, we estimate their integrals using \eqref{eq:MCin1}
\begin{equation*}
  2\int_{\Omega}|f \nabla \psi \cdot \nabla \overline{u}|
  \le 2\|\nabla \psi\|_{L^{\infty}}\|f\|_{X^{*}}\|\nabla_{x} u\|_{X}
\end{equation*}
(recall $\psi=\psi(x)$) and \eqref{eq:MCin4}
\begin{equation*}
  \int_{\Omega}|f (\Delta \psi -\phi) \overline{u}|\le
  2\||x|(\Delta \psi -\phi)\|_{L^{\infty}}
  \|f\|_{X^{*}}\|u\|_{X_{1}}.
\end{equation*}
Summing up, the integral over $\Omega$ of the right hand side
of \eqref{eq:fundest} is bounded by
\begin{equation}\label{eq:RHSint}
  C(\phi,\psi)\delta(\lambda^{+}\|u\|_{X}^{2} 
       +\|\nabla_{x} u\|_{X}^{2}+\|u\|_{X_{1}}^{2})
   +C(\phi,\psi)\delta^{-1}\|f\|_{X^{*}}^{2}
\end{equation}
with
\begin{equation}\label{eq:cphipsi}
  C(\phi,\psi)=10\|\nabla \psi\|_{L^{\infty}}+
    10\||x|(\Delta \psi -\phi)\|_{L^{\infty}}. 
\end{equation}

Consider now the left hand side of \eqref{eq:fundest}.
The term in divergence form $\nabla \cdot \Re Q_{1}$, with
\begin{equation*}
  Q_{1}=\nabla \psi|\nabla u|^{2}
  -2\nabla u(\nabla \psi \cdot 
  \nabla\overline{u})
\end{equation*}
can be handled as above by integrating first
on the cylinder $C_{R}$ and then letting $R\to+\infty$. The integral
on $S_{2}$ satisfies again \eqref{eq:liminf} and vanishes in the limit.
As to the integral on $S_{1} \subseteq \partial\Omega$, we notice that
$\nabla u$ at $\partial \Omega$ must be normal to the boundary, 
because of the Dirichlet boundary condition;
in other words, denoting the normal derivative at $\partial \Omega$ with
$\partial_{\nu}u=\nu  \cdot\nabla u$, we must have
\begin{equation*}
  \nabla u=\nu \partial_{\nu}u \ \text{\ \ at\ \  $\partial \Omega$}
\end{equation*}
so that
\begin{equation*}
  \nu \cdot\nabla Q_{1}=
  \nu \cdot \nabla \psi|\partial_{\nu}u|^{2}-
  2 \partial_{\nu}u(\nabla \psi \cdot\nu\ 
  \partial_{\nu}\overline{u})=
  -(\nu \cdot \nabla \psi)|\partial_{\nu}u|^{2}.
\end{equation*}
Thus the integral on $S_{1}$ can be written
\begin{equation*}
  I_{R}=-\int_{S_{1}}\nu \cdot \nabla \psi|\partial_{\nu}u|^{2}d \sigma
\end{equation*}
and under the second assumption \eqref{eq:asspsi} on the weight $\psi$
we obtain
\begin{equation*}
  I_{R}\ge0 \ \text{for all $R$}.
\end{equation*}
Hence we can drop $I_{R}$ from the computation, and
recalling also \eqref{eq:RHSint}
we obtain the basic integral inequality
\begin{equation}\label{eq:fundestint}
  \begin{split}
    &
    \int_{\Omega}
    \bigl[
    \nabla u(2D^{2}\psi-\phi I)\nabla \overline{u}+
    \frac12 \Delta(\phi-\Delta \psi)|u|^{2}
    +\phi \lambda|u|^{2}
    -(\nabla V \cdot \nabla \psi+\phi V)|u|^{2}\bigr]
    \le
       \\
    & \qquad \qquad
    \le
    C(\phi,\psi)\delta(\lambda^{+}\|u\|_{X}^{2} 
         +\|\nabla_{x} u\|_{X}^{2}+\|u\|_{X_{1}}^{2})
     +C(\phi,\psi)\delta^{-1}\|f\|_{X^{*}}^{2}
  \end{split}  
\end{equation}

It remains to choose the functions $\phi,\psi$ in an appropriate
way. When $\lambda>0$ we make the following choice, inspired by
 \cite{BarceloRuizVegaVilela-a}:
\begin{equation}\label{eq:psiphi1}
  \psi(x,y)=
  \begin{cases}
    |x| &\text{if $ |x|\ge R $,}\\
    \frac R2+\frac{|x|^{2}}{2R} &\text{if $ |x|< R $,}
  \end{cases}
  \qquad
  \phi(x,y)=
  \begin{cases}
    0 &\text{if $ |x|\ge R $,}\\
    \frac1R &\text{if $ |x|<R $.}
  \end{cases}
\end{equation}
Notice that assumptions \eqref{eq:assphipsi} and \eqref{eq:asspsi}
(i.e.~\eqref{eq:asspsi2})
are satisfied. We compute the quantities
relevant to our estimate: we have
\begin{equation*}
  \phi-\Delta \psi=
  \begin{cases}
    -\frac{n-1}{|x|} &\text{if $ |x|\ge R $,}\\
    -\frac{n-1}{R} &\text{if $ |x|<R $}
  \end{cases}
\end{equation*}
(with a cancelation of the singularity at $|x|=R$).
Thus we have, in distribution sense,
\begin{equation*}
  \Delta(\phi-\Delta \psi)=
  \frac{n-1}{R^{2}}\delta_{|x|=R}+
  \begin{cases}
    \frac{\mu_{n}}{|x|^{3}} &\text{if $ |x|\ge R $,}\\
    0 &\text{if $ |x|<R $,}
  \end{cases}
  \qquad
  \mu_{n}=(n-1)(n-3)
\end{equation*}
and also
\begin{equation*}
  \|\nabla \psi\|_{L^{\infty}}=1,\qquad
  \||x|(\Delta \psi-\phi)\|_{L^{\infty}}=n-1
  \qquad
  \implies
  \qquad
  C(\phi,\psi)=10n.
\end{equation*}
For the first term in \eqref{eq:fundestint} we need the 
elementary formula, valid for a radial function
$\psi=\sigma(|x|)$
\begin{equation*}
  \nabla uD^{2}\psi\nabla \overline{u}=
     \sigma''|\partial_{x} u|^{2}+
     \frac{\sigma'}{|x|}|\nabla_{x}u-\widehat{x}\ \partial_{x}u|^{2}
\end{equation*}
which implies
\begin{equation*}
  \nabla u(2D^{2}\psi-\phi I)\nabla \overline{u}=
  \begin{cases} 
    \frac2R|\nabla_{x}u-\widehat{x}
       \partial_{x}u|^{2} &\text{if $ |x|\ge R $,}\\
    \frac1R|\nabla_{x}u|^{2} &\text{if $ |x|<R $.}
  \end{cases}
\end{equation*}
Finally, the terms containing the potential $V$ are easily seen to be
positive, thanks to assumption \eqref{eq:assV}, and we can drop them.
Thus \eqref{eq:fundestint} implies
\begin{equation*}
  \begin{split}
     \frac1R\|\nabla_{x}u\|^{2}_{L^{2}(\Omega\cap\{|x|\le R\})}+
     &
     \frac{n-1}{2R^{2}}\int_{\Omega\cap\{|x|=R\}}|u|^{2}d \sigma
     +\frac \lambda R\|u\|^{2}_{L^{2}(\Omega\cap\{|x|\le R\})}\le
  \\
     & 
     \le
     10n\delta(\lambda\|u\|_{X}^{2} 
          +\|\nabla_{x} u\|_{X}^{2}+\|u\|_{X_{1}}^{2})
      +10n\delta^{-1}\|f\|_{X^{*}}^{2}
  \end{split}
\end{equation*}
and taking the sup in $R>0$ we obtain
\begin{equation*}
  \|\nabla_{x}u\|_{X}^{2}+\frac{n-1}{2}\|u\|_{X_{2}}^{2}
  +\lambda\|u\|_{X}^{2}\le
   10n\delta(\lambda\|u\|_{X}^{2} 
        +\|\nabla_{x} u\|_{X}^{2}+\|u\|_{X_{1}}^{2})
    +10n\delta^{-1}\|f\|_{X^{*}}^{2}
\end{equation*}
Recalling that the $X_{2}$ norm dominates the $X_{1}$
norm and choosing $\delta=(20n)^{-1}$
we finally obtain in the case $\lambda>0$
\begin{equation}\label{eq:final1}
  \|\nabla_{x}u\|_{X}^{2}+\|u\|_{X_{1}}^{2}
  +\lambda\|u\|_{X}^{2}\le
   400 n^{2}\|f\|_{X^{*}}^{2},\qquad
   \lambda>0.
\end{equation}

In the case $\lambda\le0$ we make a different choice of weights.
Following \cite{DanconaFanelli08-a},
we take simply $\phi \equiv0$ and we define
\begin{equation}\label{eq:defpsi}
  \psi(x)=\int_{0}^{|x|}\alpha(r)dr,\qquad
  \alpha(r)=
  \begin{cases}
    \frac1n-\frac{1}{2n(n+2)}
         \frac{R^{n-1}}{r^{n-1}} &\text{if $ r\ge R $,}\\
    \frac{1}{2n}+\frac{r}{2nR}-
         \frac{1}{2n(n+2)} \frac{r^{3}}{R^{3}} &\text{if $ r<R $.}
  \end{cases}
\end{equation}
We have now, after some elementary computations,
\begin{equation}\label{eq:lappsi}
  \Delta \psi=
  \begin{cases}
    \frac{3(n-1)}{n}\frac1r &\text{if $ r\ge R $,}\\
    \frac{1}{2R}+\frac{n-1}{nr}-
         \frac{r^{2}}{2nR^{3}} &\text{if $ r<R $,}
  \end{cases}
\end{equation}
moreover
\begin{equation*}
  \|\nabla \psi\|_{L^{\infty}}=\frac1n,\qquad
  \||x|\Delta \psi\|_{L^{\infty}}\le1-\frac1n
  \qquad \implies \qquad C(\phi,\psi)\le 10,
\end{equation*}
for $n=3$
\begin{equation*}
  -\Delta^{2}\psi=\frac{1}{R^{3}}\chi_{|x|<R}+8\pi \delta_{0}(x)
\end{equation*}
where $\delta_{0}(x)$ is the Dirac delta at 0 in the variables $x$
and $\chi_{A}$ is the characteristic function of the set $A$,
while for $n\ge4$ we have ($\mu_{n}=(n-1)(n-3)$)
\begin{equation*}
  -\Delta^{2}\psi=
  \left(
    \frac{1}{R^{3}}+\frac{\mu_{n}}{2n|x|^{3}}
  \right)\chi_{|x|<R}+
  \frac{\mu_{n}}{n|x|^{3}}\chi_{|x|\ge R}+
  \frac{n-3}{2nR^{2}}\delta_{|x|=R}
\end{equation*}
so that in all cases $n\ge3$ we have
\begin{equation*}
  -\Delta^{2}\psi\ge \frac{1}{R^{3}}\chi_{|x|<R}.
\end{equation*}
Moreover,
\begin{equation*}
  \nabla uD^{2}\psi\nabla \overline{u}\ge
  \frac{n-1}{2n(n+2)}\frac1R|\nabla_{x} u|^{2}\chi_{|x|<R}.
\end{equation*}
Thus, proceeding exactly as above, we obtain
\begin{equation*}
  \frac{n-1}{n(n+2)}
  \|\nabla_{x}u\|_{X}^{2}+\|u\|_{X_{1}}^{2}
  \le 10 \delta(\|\nabla u\|_{X}^{2}+\|u\|_{X_{1}}^{2})
  +10 \delta^{-1}\|f\|_{X^{*}}^{2}
\end{equation*}
and choosing $\delta=(40n)^{-1}$ we conclude, for $\lambda\le0$,
\begin{equation}\label{eq:final2}
  \|\nabla_{x}u\|_{X}^{2}+\|u\|_{X_{1}}^{2}\le
  800n^{2}\|f\|_{X^{*}}^{2}.
\end{equation}

We collect \eqref{eq:final1} and \eqref{eq:final2} in the estimate,
valid for all $\lambda\in \mathbb{R}$,
\begin{equation}\label{eq:final3}
  \|\nabla_{x}u\|_{X}^{2}+\|u\|_{X_{1}}^{2}
  +\lambda^{+}\|u\|_{X}^{2}\le
   800 n^{2}\|f\|_{X^{*}}^{2}.
\end{equation}
As a last step, we show that the factor $\lambda^{+}$ in \eqref{eq:final3}
can be improved to $|\lambda|+|\epsilon|$. First of all, recalling
\eqref{eq:basic2}, and using \eqref{eq:final3}, we see that
\begin{equation}\label{eq:final4}
  |\epsilon|\|u\|^{2}_{X}\le 4(\|f\|^{2}_{X^{*}}+\|\nabla u\|_{X})
  \| u\|_{X_{1}}\le 3320 n^{2}\|f\|^{2}_{X^{*}}.
\end{equation}

Assume now $\lambda=-\lambda^{-}\le0$.
We multiply the resolvent equation \eqref{eq:reseq}
by $\overline{u}$ and take real parts, obtaining
\begin{equation*}
  |\nabla u|^{2}+\lambda^{-}|u|^{2}+V|u|^{2}=
  \Re(f \overline{u})+\frac12 \Delta|u|^{2};
\end{equation*}
then we multiply by a weight function $\mu(x)$ and we get
\begin{equation*}
  \mu |\nabla u|^{2}+(\lambda^{-}+V)\mu|u|^{2}=
  \Re(\mu f \overline{u})+\frac12 \Delta\mu |u|^{2}+
  \nabla \cdot(2^{-1}\nabla(\mu|u|^{2})).
\end{equation*}
We now integrate on $\Omega$ as above; the term in divergence form
vanishes by the Dirichlet b.c., and we obtain, using the positivity
of $V$,
\begin{equation}\label{eq:lastid}
  \int_{\Omega}\mu(|\nabla u|^{2}+
  \lambda^{-}|u|^{2})\le 
  \int_{\Omega}\mu|f \overline{u}|+
  \frac12\int_{\Omega}\Delta\mu|u|^{2}.
\end{equation}
We now choose $\mu=\Delta\psi$ with $\psi$ defined as in \eqref{eq:defpsi}.
Notice that $\Delta\mu=\Delta^{2}\psi\le0$ so we can drop
the last term from the computation; on the other hand
\begin{equation*}
  \Delta\psi\ge \frac1{2R}\chi_{|x|<R},\qquad
  \||x|\Delta\psi\|_{L^{\infty}}\le1
\end{equation*}
and recalling property \eqref{eq:MCin4} we obtain
\begin{equation*}
  \frac{1}{2R}\int_{|x|<R}(|\nabla u|^{2}+\lambda^{-}|u|^{2})\le 
  2\|f\|_{X^{*}}\|u\|_{X_{1}}.
\end{equation*}
Taking the sup in $R>0$ this gives
\begin{equation}\label{eq:final5}
  \|\nabla u\|_{X}^{2}+\lambda^{-}\|u\|_{X}^{2}\le 
  4\|f\|_{X^{*}}\|u\|_{X_{1}}\le 120 n\|f\|_{X^{*}}^{2}
\end{equation}
again by \eqref{eq:final3}.

Collecting \eqref{eq:final5}, \eqref{eq:final4} and \eqref{eq:final3}
we conclude the proof of \eqref{eq:fundresest}.

\end{proof}

\begin{remark}\label{rem:MCweight}
  When $z=\lambda+i \epsilon$ does not belong to the spectrum of
  the selfadjoint operator
  $H=-\Delta+V$ with Dirichlet b.c.~on $L^{2}(\Omega)$
  (this includes some cases when $\epsilon=0$),
  given an $f\in L^{2}(\Omega)$, we can represent
  the solution of \eqref{eq:reseq} as $u=R(z)f$, where
  $R(z)=(H-z)^{-1}$. Since we know
  that $u\in H^{1}_{0}(\Omega)$, all the preceding computations
  apply and in particular estimate \eqref{eq:fundresest} holds.
  As a consequence, using \eqref{eq:MCtoweight} and 
  \eqref{eq:weighttoMC}, we can write for all $R,S>0$,
  \begin{equation}\label{eq:simpleresest}
    \|\xxR^{-3}R(z)f\|_{L^{2}(\Omega)}\le
    2^{13}n\|\xxS f\|_{L^{2}(\Omega)}.
  \end{equation}
  Thus \eqref{eq:simpleresest}
  is in fact a weighted $L^{2}$ estimate for the resolvent 
  $R(z)$. By duality we have the equivalent estimate
  \begin{equation}\label{eq:simpledual}
    \|\xxR^{-1}R(z)f\|_{L^{2}(\Omega)}\le
    c_{n}\|\xxS^{3} f\|_{L^{2}(\Omega)}
  \end{equation}
  and by (complex) interpolation we have also
  \begin{equation*}
    \|\xxR^{-2}R(z)f\|_{L^{2}(\Omega)}\le
    c_{n}\|\xxS^{2} f\|_{L^{2}(\Omega)},
  \end{equation*}
  uniformly in $z\not\in \sigma(H)$, which we shall
  write more symmetrically as follows:
  \begin{equation}\label{eq:simpleinterp}
    \|\xxR^{-2}R(z)\xxS^{-2}f\|_{L^{2}(\Omega)}\le
    c_{n}\|f\|_{L^{2}(\Omega)}.
  \end{equation}
  A similar computation, using the other two terms in 
  \eqref{eq:fundresest}, shows that
  \begin{equation}\label{eq:gradest}
    \|\xxR^{-1}\nabla_{x} R(z)\xxS^{-1}f\|_{L^{2}(\Omega)}
    +|z|^{1/2}
    \|\xxR^{-1} R(z)\xxS^{-1}f\|_{L^{2}(\Omega)}
    \le
    c_{n}\|f\|_{L^{2}(\Omega)}
  \end{equation}
  uniformly in $z\not\in \sigma(H)$. In particular this applies to
  $z=-\delta$ for all $\delta>0$ since the operator $H$ is positive.
  
  At this point we need the following elementary
  
\begin{lemma}\label{lem:weights}
  If a linear operator $A$ satisfies for all $R,S>0$ the estimate
  \begin{equation}\label{eq:basest}
    \|\xxR^{-\gamma}A\xxS^{-\gamma}u\|_{L^{2}}
    \le C_{0}\|u\|_{L^{2}}
  \end{equation}
  with a constant independent of $R,S$, then it satisfies also,
  for all $\epsilon>0$, the estimate
  \begin{equation}\label{eq:basest1}
    \|\xx^{-\frac{\gamma}{2}-\epsilon}A
       \xx^{-\frac{\gamma}{2}-\epsilon}u\|_{L^{2}}\le 
       C_{0}C(\gamma,\epsilon)\|u\|_{L^{2}}.
  \end{equation}
\end{lemma}

\begin{proof}
  Write \eqref{eq:basest} in the form
  \begin{equation}\label{eq:basest2}
    \|\xxR^{-\gamma}Av\|_{L^{2}}
    \le C_{0}\|\xxS^{\gamma}v\|_{L^{2}},
  \end{equation}
  decompose $v=v_{0}+\sum_{j\ge1} v_{j}$, with $v_{j}$ supported in 
  $2^{j-1}\le|x|<2^{j}$ for $j\ge1$ and $v_{0}$ in $|x|<1$,
  apply the \eqref{eq:basest2} to each $v_{j}$
  with $S=2^{j}$, and sum over $j$ (all norms in the rest of the proof
  are in $L^{2}$):
  \begin{equation*}
    \|\xxR^{-\gamma}Av\|
    \le \|\xxR^{-\gamma}Av_{0}\|+
    \sum_{j}\|\xxR^{-\gamma}Av_{j}\|
    \le C_{0}\|\xx^{\gamma}v_{0}\|+
    C_{0}\sum\|\xx_{2^{j}}^{\gamma}v_{j}\|.
  \end{equation*}
  Now notice that for $j\ge1$
  \begin{equation*}
    \xx_{2^{j}}^{2\gamma}
    =\left(2^{j}+\frac{|x|^{2}}{2^{j}}\right)^{\gamma}
    \le 2^{\gamma}2^{\gamma j}\le
    2^{2 \gamma}|x|^{\gamma}\le
    2^{2 (\gamma+\epsilon)}2^{-2 \epsilon j}|x|^{\gamma+2 \epsilon}
    \quad\text{on the support of $v_{j}$}
  \end{equation*}
  so that
  \begin{equation*}
    \|\xxR^{-\gamma}Av\|\le
    C_{0}2^{\gamma}
    \|v_{0}\|+C_{0}2^{\gamma+\epsilon}
    \sum_{j\ge1} 2^{-\epsilon j}
    \||x|^{\frac{\gamma}{2}+\epsilon}v_{j}\|\le
    C_{0}C(\gamma,\epsilon)\|\xx^{\frac{\gamma}{2}+\epsilon}v\|
  \end{equation*}
  by Cauchy-Schwarz. Using \eqref{eq:weight1} we obtain
  \eqref{eq:basest1}.
\end{proof}
  
  In particular, applying the Lemma to \eqref{eq:simpleinterp}
  and to \eqref{eq:gradest} we obtain the estimates,
  valid for all $\epsilon>0$:
  \begin{equation}\label{eq:simpleinterpbis}
    \|\xx^{-1-\epsilon}R(z)\xx^{-1-\epsilon}f\|_{L^{2}(\Omega)}\le
    c_{n,\epsilon}\|f\|_{L^{2}(\Omega)},
  \end{equation}
  \begin{equation}\label{eq:gradestbis}
    \|\xx^{-\frac12-\epsilon}\nabla_{x} 
         R(z)\xx^{-\frac12-\epsilon}f\|_{L^{2}(\Omega)}\le
    c_{n,\epsilon}\|f\|_{L^{2}(\Omega)},
  \end{equation}
  \begin{equation}\label{eq:gradestter}
    |z|^{1/2}
    \|\xx^{-\frac12-\epsilon}
         R(z)\xx^{-\frac12-\epsilon}f\|_{L^{2}(\Omega)}\le
    c_{n,\epsilon}\|f\|_{L^{2}(\Omega)}.
  \end{equation} 
\end{remark}


\section{Smoothing estimates}\label{sec:smoothing}  

The concept of \emph{$H$-smoothing} was introduced by
Kato \cite{Kato65-a} in the context of scattering theory, and
its usefulness for dispersive equations was
revealed in \cite{RodnianskiSchlag04-a}. An operator $A$
is \emph{$H$-smooth} (actually, supersmooth) whenever one
of the two equivalent estimates \eqref{eq:resH}, \eqref{eq:smoH1}
in the following theorem holds. We shall use a version of the
result adapted to the applications we have in mind; for a more
complete reference see \cite{ReedSimon75-a},
\cite{Mochizuki-preprint}

\begin{theorem}[Kato]\label{the:katosmoo}
  Assume $K$ is a selfadjoint operator in a Hilbert space
  $\mathcal{H}$, let $\mathcal{R}(z)=(K-z)^{-1}$ 
  be its resolvent operator
  for $z\in \mathbb{C}\setminus \mathbb{R}$, and let $A$ be
  a densely defined closed operator from $\mathcal{H}$ to a
  second Hilbert space $\mathcal{H}_{1}$ with $D(A)\supseteq D(K)$.
  
  Assume that $A,\mathcal{R}(z)$ satisfy the estimate
  \begin{equation}\label{eq:resH}
    \sup_{z\not\in \mathbb{R}}\|A \mathcal{R}(z)A^{*}f\|_{\mathcal{H}_{1}}
    \le c_{0}^{2}\|f\|_{\mathcal{H}_{1}}
  \end{equation}
  for all $f\in D(A^{*})$. Then the following estimates hold:
  \begin{equation}\label{eq:smoH1}
    \|Ae^{itK}f\|_{L^{2}_{t}\mathcal{H}_{1}}\le
    c_{0}\|f\|_{\mathcal{H}},
  \end{equation}
  \begin{equation}\label{eq:smooH2}
    \left\|
      \int_{0}^{t} A e^{i(t-s)K}A^{*}h(s)ds
    \right\|_{L^{2}_{t}\mathcal{H}_{1}}
    \le c_{0}^{2}\|h\|_{L^{2}_{t}\mathcal{H}_{1}}
  \end{equation}
  for all $f\in \mathcal{H}$, $h\in L^{2}_{t}\mathcal{H}_{1}$.
  
  Estimate \eqref{eq:smoH1} still holds when
  \eqref{eq:resH} is replaced by the weaker assumption
  \begin{equation}\label{eq:resH2}
    \sup_{z\not\in \mathbb{R}}\|A\Im(
       \mathcal{R}(z))A^{*}
    f\|_{\mathcal{H}_{1}}
    \le c_{0}^{2}\|f\|_{\mathcal{H}_{1}},
  \end{equation}
  where we use the notation $\Im T=(2i)^{-1}(T-T^{*})$.
\end{theorem}

Recalling \eqref{eq:simpleinterp} in Remark \ref{rem:MCweight},
we see that with the choices
\begin{equation*}
  \mathcal{H}=\mathcal{H}_{1}=
  L^{2}(\Omega),\qquad
  K=H=-\Delta+V(x,y),\qquad
  A=\xx^{-1-\epsilon}
\end{equation*}
estimate \eqref{eq:simpleinterpbis} reduces precisely to \eqref{eq:resH}.
Thus from Theorem \ref{the:katosmoo} and \eqref{eq:simpleinterpbis}
we obtain immediately the  following smoothing estimates for the
Schr\"odinger flow associated to the operator $H=-\Delta+V(x,y)$:

\begin{theorem}\label{the:smoosch}
  Let the domain
  $\Omega \subseteq \mathbb{R}^n_{x}\times\mathbb{R}^{m}_{y}$, 
  $n\ge3$, $m\ge1$ be
  repulsive with respect to the $x$ variables,
  with a Lipschitz boundary.
  Assume the operator
  $H=-\Delta+V(x,y)$ with Dirichlet b.c. is selfadjoint
  on $L^{2}(\Omega)$. Finally,
  assume that the potential $V$
  satisfies on $\Omega$ the inequalities
  \begin{equation}\label{eq:assVdel}
    V(x,y)\ge0,
    \qquad
    -\partial_{x}(|x|V(x,y))\ge0.
  \end{equation}
  Then the Schr\"odinger flow associated to $H$ satisfies
  the following smoothing estimates: for any $\epsilon>0$,
  \begin{equation}\label{eq:smoHV1}
    \|\xx^{-1-\epsilon}e^{itH}f\|_{L^{2}_{t}L^{2}(\Omega)}\lesssim
    \|f\|_{L^{2}(\Omega)},
  \end{equation}
  \begin{equation}\label{eq:smooHV2}
    \left\|\xx^{-1-\epsilon}
      \int_{0}^{t} e^{i(t-s)H}F(s)ds
    \right\|_{L^{2}_{t}L^{2}(\Omega)}
    \lesssim \|\xx^{1+\epsilon} F\|_{L^{2}_{t}L^{2}(\Omega)}
  \end{equation}
  for all $f(x,y)\in L^{2}(\Omega)$ and $F(t,x,y)$ with
  $\xx^{1+\epsilon} F\in L^{2}_{t}L^{2}(\Omega)$.
\end{theorem}

We can obtain an estimate also for the derivatives of
$e^{itH}f$, with a gain of a half derivative, by a different
choice of the operator $A$ and some functional analytic arguments;
to this end we must introduce suitable functional spaces.

For functions on $\mathbb{R}^{n+m}$ and
$z\in \mathbb{C}$, we introduce the operators acting only on the
$x$ variables
\begin{equation*}
  |D_{x}|^{z}f(x,y)=(2\pi)^{-n}
     \int_{\mathbb{R}^{n}}|\xi|^{z}
     \widehat{f}(\xi,y)e^{i\xi x}dx,
\end{equation*}
\begin{equation*}
  \langle D_{x}\rangle ^{z}f(x,y)=(2\pi)^{-n}
     \int_{\mathbb{R}^{n}}\xxi^{z}
     \widehat{f}(\xi,y)e^{i\xi x}dx,
\end{equation*}
where $\widehat{f}(\xi,y)$ is the Fourier transform of $f(x,y)$
with respect to the variable $x$ only.
By standard calculus we have the equivalence
\begin{equation*}
  \||D_{x}|f\|_{L^{2}(\mathbb{R}^{n+m})}\simeq
  \|\nabla_{x} f\|_{L^{2}(\mathbb{R}^{n+m})}.
\end{equation*}
We introduce the norms, and the corresponding Hilbert spaces,
\begin{equation}\label{eq:xsobolev}
  \|f\|_{\dot H^{s,0}}=
  \||D_{x}|^{s}f\|_{L^{2}(\mathbb{R}^{n+m})},\qquad
  \|f\|_{H^{s,0}}=
  \|\langle D_{x}\rangle^{s}f\|_{L^{2}(\mathbb{R}^{n+m})}.
\end{equation}
Notice that,
if the boundary of $\Omega$ satisfies a uniform Lipschitz condition,
the extension as 0 of a function $f\in H^{1}_{0}(\Omega)$
to all of $\mathbb{R}^{n+m}$ gives a function
$Ef\in H^{1}\mathbb{R}^{n+m}$ with the same norm.
Thus for $f\in H^{1}_{0}(\Omega)$ and $0\le\Re z\le1$
we can extend the definition of the operators as
\begin{equation*}
  |D_{x}|^{z}f=|D_{x}|^{z}Ef,\qquad
  \langle D_{x}\rangle^{z}f=|D_{x}|^{z}Ef.
\end{equation*}
By density of $C^{\infty}_{c}(\Omega)$ in $H^{1}_{0}(\Omega)$ we
obtain also that
\begin{equation}\label{eq:equiv}
  \||D_{x}|f\|_{L^{2}(\Omega)}\simeq
  \|\nabla_{x} f\|_{L^{2}(\Omega))}.
\end{equation}
Recall now the estimate ($y\in \mathbb{R}$)
\begin{equation}\label{eq:riesz}
  \|\xx^{-s}|D_{x}|^{1+iy}f\|_{L^{2}(\mathbb{R}^{n+m})}\simeq
  \|\xx^{-s}\nabla_{x} f\|_{L^{2}(\mathbb{R}^{n+m})},
  \qquad
  s>-\frac n2
\end{equation}
which holds since the Riesz operators  $\partial_{x_{j}}|D_{x}|^{-1}$
and the operators $|D_{x}|^{iy}$
are bounded in weighted $L^{2}$ with $A_{2}$ weights, and
$\xx^{-s}\in A_{2}(\mathbb{R}^{n})$ for $s>-n/2$;
notice also that the constant in the estimate depends on
$y\in \mathbb{R}$ but with a polynomial growth as $|y|\to \infty$
(see \cite{Stein93-a} for the general theory of singular integrals
in weighted $L^{2}$ spaces, and more specifically
\cite{SikoraWright01-a}, \cite{CacciafestaDAncona09-a} 
for the polynomial growth of the norms). The estimate extends to
\begin{equation}\label{eq:equivw}
  \|\xx^{-s}|D_{x}|^{1+iy}f\|_{L^{2}(\Omega)}\simeq
  \|\xx^{-s}\nabla_{x} f\|_{L^{2}(\Omega))},\qquad
    s>-\frac n2,\quad
    f\in H^{1}_{0}(\Omega)
\end{equation}
by a density argument as above.

As a consequence of \eqref{eq:equivw}, 
\eqref{eq:gradestbis} implies the estimate
\begin{equation}\label{eq:gradest2}
  \|\xx^{-1/2-\epsilon}|D_{x}|^{1+iy}
   R(z)\xx^{-1/2-\epsilon}f\|_{L^{2}(\Omega)}
  \le
  c_{n,\epsilon}\|f\|_{L^{2}(\Omega)}
\end{equation}
which by duality is equivalent to
\begin{equation}\label{eq:gradest2dual}
  \|\xx^{-1/2-\epsilon}
   R(z)|D_{x}|^{1+iy}\xx^{-1/2-\epsilon}f\|_{L^{2}(\Omega)}
  \le
  c_{n,\epsilon}\|f\|_{L^{2}(\Omega)}
\end{equation}
Thus by complex interpolation 
for the analytic family of operators
$T_{z}=$
we also obtain the estimate
\begin{equation}\label{eq:gradesthalf}
  \|\xx^{-1/2-\epsilon}|D_{x}|^{1/2}
   R(z)|D_{x}|^{1/2}\xx^{-1/2-\epsilon}f\|_{L^{2}(\Omega)}
  \le
  c_{n,\epsilon}\|f\|_{L^{2}(\Omega)}
\end{equation}

Now we make the following choice:
\begin{equation*}
  \mathcal{H}=\dot H^{1/2,0}(\Omega),\qquad
  \mathcal{H}_{1}=L^{2}(\Omega),\qquad
  H=-\Delta+V(x,y)
\end{equation*}
where the space $\dot H^{1/2,0}(\Omega)$ is defined as the
completion of $C^{\infty}_{c}(\Omega)$ in the norm
\begin{equation*}
  \|f\|_{\dot H^{1/2,0}(\Omega)}=
  \||D_{x}|^{1/2}f\|_{L^{2}(\Omega)},
\end{equation*}
The closed unbounded operator $A:\mathcal{H}\to \mathcal{H}_{1}$
is now defined as
\begin{equation*}
  A=\xx^{-1/2-\epsilon}|D_{x}|
\end{equation*}
and its adjoint $A^{*}$ is computed as follows
\begin{equation*}
  \begin{split}
  (Af,g)_{\mathcal{H}_{1}}=
  &
  (\xx^{-1/2-\epsilon}|D_{x}|f,g)_{L^{2}(\Omega)}=
  (|D_{x}|f,\xx^{-1/2-\epsilon}g)_{L^{2}(\Omega)}=
  \\
  &
  =(|D_{x}|^{1/2}f,|D_{x}|^{1/2}\xx^{-1/2-\epsilon}g)_{L^{2}(\Omega)}=
  (f,\xx^{-1/2-\epsilon}g)_{\mathcal{H}}=
  (f,A^{*}g)_{\mathcal{H}}.
  \end{split}
\end{equation*}
With these choices, estimate \eqref{eq:gradestbis} takes precisely 
the form \eqref{eq:resH} and Kato theory applies. We obtain
the following
\begin{theorem}\label{the:smoosch2}
  Assume $\Omega$, $V$, $H$ as in Theorem \ref{the:smoosch}.
  Then the Schr\"odinger flow associated to $H$ satisfies
  the smoothing estimates
  \begin{equation}\label{eq:smoHV3}
    \|\xx^{-1/2-\epsilon}\nabla_{x} e^{itH}f\|_{L^{2}_{t}L^{2}(\Omega)}
    \lesssim
    \||D_{x}|^{1/2}f\|_{L^{2}(\Omega)},
  \end{equation}
  \begin{equation}\label{eq:smooHV4}
    \left\|\xx^{-1/2-\epsilon}
      \int_{0}^{t}  \nabla_{x} e^{i(t-s)H}F(s)ds
    \right\|_{L^{2}_{t}L^{2}(\Omega)}
    \lesssim \|\xx^{1/2+\epsilon} F\|_{L^{2}_{t}L^{2}(\Omega)}
  \end{equation}
  for all $f(x,y)\in H^{1}_{0}(\Omega)$ and $F(t,x,y)$ with
  $\xx^{1/2+\epsilon} F\in L^{2}_{t}L^{2}(\Omega)$.
\end{theorem}

Notice that a different choice is possible: namely, if we set
\begin{equation*}
  \mathcal{H}=
  \mathcal{H}_{1}=L^{2}(\Omega),\qquad
  H=-\Delta+V(x,y)
\end{equation*}
and
\begin{equation*}
  A=\xx^{-1/2-\epsilon}|D_{x}|^{1/2},\qquad 
  A^{*}=|D_{x}|^{1/2}\xx^{-1/2-\epsilon}
\end{equation*}
we obtain the (essentially equivalent) result:
\begin{theorem}\label{the:smoVH3bis}
  Assume $\Omega$, $V$, $H$ as in Theorem \ref{the:smoosch}.
  Then the Schr\"odinger flow associated to $H$ satisfies
  the smoothing estimates
  \begin{equation}\label{eq:smoHV3bis}
    \|\xx^{-1/2-\epsilon}|D_{x}|^{1/2} e^{itH}f\|_{L^{2}_{t}L^{2}(\Omega)}
    \lesssim
    \|f\|_{L^{2}(\Omega)},
  \end{equation}
  \begin{equation}\label{eq:smooHV4bis}
    \left\|\xx^{-1/2-\epsilon}
      \int_{0}^{t} |D_{x}|^{1/2} e^{i(t-s)H}F(s)ds
    \right\|_{L^{2}_{t}L^{2}(\Omega)}
    \lesssim \|\xx^{1/2+\epsilon} |D_{x}|^{-1/2}F\|_{L^{2}_{t}L^{2}(\Omega)}
  \end{equation}
  for all $f(x,y)\in L^{2}(\Omega)$ and $F(t,x,y)$ with
  $\xx^{1/2+\epsilon} F\in L^{2}_{t}L^{2}(\Omega)$.
\end{theorem}

Handling the wave and Klein-Gordon equations requires some additional
effort. We start from the standard representation
\begin{equation}\label{eq:matrix1}
  K=
  \begin{pmatrix}
   0 & 1 \\
   H &  0
  \end{pmatrix}
  \qquad\implies \qquad
  \exp(itK)=
    \begin{pmatrix}
     \cos(t \sqrt{H}) & \frac{i}{\sqrt{H}}\sin(t \sqrt{H}) \\
      i \sqrt{H}\sin(t \sqrt{H})&  \cos(t \sqrt{H})
    \end{pmatrix}
\end{equation}
so that
\begin{equation}\label{eq:matrix2}
  e^{itK}
  \begin{pmatrix}
    f \\
    \sqrt{H}f  
  \end{pmatrix}=
  \begin{pmatrix}
   e^{it \sqrt{H}}f  \\
   \sqrt{H}  e^{it \sqrt{H}}f 
  \end{pmatrix}
\end{equation}
is the flow associated to the wave equation
\begin{equation*}
  u_{tt}+Hu=0.
\end{equation*}
Now we choose
\begin{equation*}
  \mathcal{H}=D(\sqrt{H})\times L^{2}(\Omega),\qquad
  \mathcal{H}_{1}=
  L^{2}(\Omega),\qquad
  H=-\Delta+V(x,y)
\end{equation*}
with $K$ as in \eqref{eq:matrix1},
and $A:\mathcal{H}\to L^{2}(\Omega)$
defined by
\begin{equation*}
  A
  \begin{pmatrix}
   f \\
   g  
  \end{pmatrix}
  =\xx^{-1/2-\epsilon}H^{1/2}f
  \qquad \implies \qquad
  A^{*}g=
  \begin{pmatrix}
   H^{-1/2}\xx^{-1/2-\epsilon}g \\
    0
  \end{pmatrix}.
\end{equation*}
Then the resolvent $\mathcal{R}(z)=(K-z)^{-1}$ 
can be written in terms of the resolvent $R(z)=(H-z)^{-1}$
as
\begin{equation}\label{eq:resres}
  \mathcal{R}(z)=
  \begin{pmatrix}
   zR(z^{2}) & R(z^{2}) \\
   HR(z^{2}) &  zR(z^{2})
  \end{pmatrix}.
\end{equation}
Thus we see that, in order to apply the Kato theory to $e^{itK}$,
we need to prove that the following operator is bounded
on $L^{2}(\Omega)$, uniformly in $z\not\in \mathbb{R}$:
\begin{equation*}
  Q(z)=
  A \mathcal{R}(z)A^{*}\equiv
  \xx^{-1/2-\epsilon} zR(z^{2})\xx^{-1/2-\epsilon}.
\end{equation*}
This is precisely what is expressed by estimate \eqref{eq:gradestter}.
Then by Theorem
\ref{the:katosmoo} we obtain
\begin{equation*}
  \left\|
  A
  e^{itK}
  \begin{pmatrix}
    f \\
    \sqrt{H}f  
  \end{pmatrix}
  \right\|_{L^{2}_{t}\mathcal{H}_{1}}\lesssim
  \left\|
  \begin{pmatrix}
    f \\
    \sqrt{H}f  
  \end{pmatrix}
  \right\|_{\mathcal{H}}
\end{equation*}
which means
\begin{equation*}
  \|\xx ^{-1/2-\epsilon}H^{1/2}e^{it \sqrt{H}}f\|_{L^{2}_{t}L^{2}(\Omega)}
  \lesssim \|H^{1/2}f\|_{L^{2}(\Omega)}
\end{equation*}
or equivalently
\begin{equation*}
  \|\xx ^{-1/2-\epsilon}e^{it \sqrt{H}}f\|_{L^{2}_{t}L^{2}(\Omega)}
  \lesssim \|f\|_{L^{2}(\Omega)}.
\end{equation*}
A similar estimate is obtained for the Duhamel term.
All the previous computations are valid if
we replace the operator $H$ with $H+\mu^{2}$, $\mu\ge0$;
this gives an analogous estimate for the flow 
$e^{it \sqrt{H+\mu^{2}}}$ associated
to the Klein-Gordon equation.
In conclusion, we have proved:

\begin{theorem}\label{the:smoowave}
  Let $\mu\ge0$ and
  assume $\Omega$, $V$, $H$ as in Theorem \ref{the:smoosch}.
  Then the wave flow associated to $H+\mu^{2}$ satisfies
  the smoothing estimates
  \begin{equation}\label{eq:smowaveHV1}
    \|\xx^{-1/2-\epsilon}
    e^{it \sqrt{H+\mu^{2}}}f\|_{L^{2}_{t}L^{2}(\Omega)}
    \lesssim
    \|f\|_{L^{2}(\Omega)},
  \end{equation}
  \begin{equation}\label{eq:smoowaveHV2}
    \left\|\xx^{-1/2-\epsilon}
      \int_{0}^{t} e^{i(t-s)\sqrt{H+\mu^{2}}}F(s)ds
    \right\|_{L^{2}_{t}L^{2}(\Omega)}
    \lesssim \|\xx^{1/2+\epsilon} F\|_{L^{2}_{t}L^{2}(\Omega)}
  \end{equation}
  for all $f(x,y)\in L^{2}(\Omega)$ and $F(t,x,y)$ with
  $\xx^{1/2+\epsilon} F\in L^{2}_{t}L^{2}(\Omega)$.
\end{theorem}


\section{Strichartz estimates for the Schr\"odinger equation}\label{sec:strichartz_estimates}  

From now on we reduce to the simpler situation when
the domain $\Omega$, besides being $x$-repulsive, is a compactly
supported perturbation of a product domain. More precisely
we assume that there exist a constant
$M$ and an open set $\omega \subseteq \mathbb{R}^{m}$ such that
\begin{equation}\label{eq:assO}
  \Omega\cap \left\{(x,y):|x|>M\right\}=
  (\mathbb{R}^{n}\times\omega)\cap \left\{(x,y):|x|>M\right\}.
\end{equation}

We recall the estimates proved in Example \ref{exa:schroflat}
in the flat case
\begin{equation}\label{eq:strichflat}
  \|e^{it \Delta}f\|_{L^{2}_{t}L^{2}_{y}L^{\frac{2n}{n-2}}_{x}}
  \lesssim \|f\|_{L^{2}(\Omega)},\qquad
  \left\|
    \int_{0}^{t}e^{i(t-s)\Delta}F(s)ds
  \right\|_{L^{2}_{t}L^{2}_{y}L^{\frac{2n}{n-2}}_{x}}
  \lesssim
  \|F\|_{L^{2}_{t}L^{2}_{y}L^{\frac{2n}{n+2}}_{x}}
\end{equation}
where $\Delta$ is the Dirichlet Laplacian on 
$\mathbb{R}^{n}\times \omega$.
In the following, we shall also need a mixed Strichartz-smoothing
nonhomogeneous estimate, which follows like \eqref{eq:strichflat}
from a corresponding estimate on the whole space. Indeed, Ionescu
and Kenig proved that for the standard Laplace operator on 
$\mathbb{R}^{n}$, $n\ge3$, one has
\begin{equation}\label{eq:IK}
  \left\|\int_{0}^{t}e^{i(t-s)\Delta}F(s)ds\right\|
    _{L^{2}_{t}L^{\frac{2n}{n+2}}_{x}}\lesssim
  \|\xx^{1/2+\epsilon}|D|^{-1/2}F\|_{L^{2}_{t}L^{2}_{x}}
\end{equation}
(see Lemma 3 in \cite{IonescuKenig05-a}, which is actually the dual
form of \eqref{eq:IK}, and in a sharper version). By mimicking the
proof of \eqref{eq:strichflat} we obtain the following
mixed estimate on a flat waveguide:
\begin{equation}\label{eq:IKflat}
  \left\|\int_{0}^{t}e^{i(t-s)\Delta}F(s)ds\right\|
    _{L^{2}_{t}L^{2}_{y}L^{\frac{2n}{n+2}}_{x}}\lesssim
  \|\xx^{1/2+\epsilon}|D_{x}|^{-1/2}F\|_{L^{2}_{t}L^{2}_{y}L^{2}_{x}}
\end{equation}
where again $\Delta$ denotes the Dirichlet Laplacian on 
$\mathbb{R}^{n}\times \omega$.

Assume now the domain $\Omega$ is repulsive
with respect to $x$ and satisfies in addition the condition 
\eqref{eq:assO}, and let $u(t,x,y)$ be a solution 
on $\Omega$ of the equation 
\begin{equation}\label{eq:schreqV}
  iu_{t}-\Delta u+V(x,y)u=0,\qquad
  u(0,x,y)=f(x,y)
\end{equation}
Recall that by \eqref{eq:smoHV1}, \eqref{eq:smoHV3} 
and \eqref{eq:smoHV3bis} $u$ satisfies
\begin{equation}\label{eq:smooforu}
  \|\xx^{-1-\epsilon}u\|_{L^{2}_{t}L^{2}(\Omega)}
  \lesssim\|f\|_{L^{2}(\Omega)}
  ,\qquad
  \|\xx^{-1/2-\epsilon}\nabla u\|_{L^{2}_{t}L^{2}(\Omega)}
  \lesssim\||D_{x}|^{1/2}f\|_{L^{2}(\Omega)}
\end{equation}
and
\begin{equation}\label{eq:smooforubis}
  \|\xx^{-1/2-\epsilon}|D_{x}|^{1/2} u\|_{L^{2}_{t}L^{2}(\Omega)}
  \lesssim\|f\|_{L^{2}(\Omega)}.
\end{equation}
Fix a cutoff function $\chi(x)$ equal to 1 on the ball
$B(0.M)$ and vanishing outside $B(0,M+1)$ and split the solution as
\begin{equation*}
  u=v+w,\qquad v(t,x,y)=\chi(x)u(t,x,y),\quad
  w(t,x,y)=(1-\chi(x))u(t,x,y).
\end{equation*}
Then $w$ is a solution of the following
Schr\"odinger equation
\begin{equation}\label{eq:perturb1}
  iw_{t}-\Delta w=G_{1}+G_{2},\qquad
  G_{1}=-V(x,y)(1-\chi)u+\Delta_{x}\chi \ u,\quad
  G_{2}=2 \nabla_{x} \chi \cdot \nabla_{x}u,
\end{equation}
\begin{equation*}
  w(0,x,y)=(1-\chi(x))f(x,y)
\end{equation*}
on $\mathbb{R}^{n}\times \omega$ with Dirichlet boundary conditions.
We can now represent $w(t,x,y)$ as
\begin{equation*}
  w=e^{it \Delta}(1-\chi)f+
   i\int_{0}^{t}e^{i(t-s)\Delta}G_{1}(s)ds+
   i\int_{0}^{t}e^{i(t-s)\Delta}G_{2}(s)ds
   \equiv I+II+III.
\end{equation*}
We plan to use estimates \eqref{eq:strichflat} on the first two terms
and \eqref{eq:IKflat} on the third one.
The $L^{2}_{t}L^{2}_{y}L^{\frac{2n}{n-2}}_{x}$ norm of the
first term $I$ is estimated directly using \eqref{eq:strichflat}.
Again by \eqref{eq:strichflat}, the 
$L^{2}_{t}L^{2}_{y}L^{\frac{2n}{n-2}}_{x}$ norm of $II$ is estimated
using H\"older's inequality as follows
\begin{equation}\label{eq:G1}
  \|\Delta_{x}\chi(x)u\|_{L^{2}_{t}L^{2}_{y}L^{\frac{2n}{n+2}}_{x}}
  \lesssim
  \|\xx^{-1-\epsilon} u\|_{L^{2}_{t}L^{2}(\Omega)}
  \lesssim\|f\|_{L^{2}(\Omega)},
\end{equation}
and
\begin{equation}\label{eq:G3}
  \|Vu\|_{L^{2}_{t}L^{2}_{y}L^{\frac{2n}{n+2}}_{x}}\le
  \|\xx^{1+\epsilon} V\|_{L^{2}_{y}L^{n}_{x}}
  \|\xx^{-1-\epsilon}u\|_{L^{2}_{t}L^{2}(\Omega)}
  \lesssim
  \|\xx^{1+\epsilon} V\|_{L^{2}_{y}L^{n}_{x}}
  \|f\|_{L^{2}(\Omega)}.
\end{equation}
using the smoothing estimate \eqref{eq:smoHV1} in both cases.
For the third term $III$, on the other hand, we use the
mixed estimate \eqref{eq:IKflat} so that
\begin{equation*}
  \|III\|_{L^{2}_{t}L^{2}_{y}L^{\frac{2n}{n-2}}_{x}}\lesssim
  \|\xx^{1/2+\epsilon}|D_{x}|^{-1/2}
  (\nabla_{x}\chi \cdot \nabla_{x}u)\|_{L^{2}_{t}L^{2}_{y}L^{2}_{x}}.
\end{equation*}
Let now $\psi(x)$ be a cutoff function supported in $|x|\le M+3$
and equal to 1 on $|x|\le M+1$
(note $\chi$ is supported in $B(0,M+1)$) and recall the explicit
formula
\begin{equation*}
  |D_{x}|^{-1/2}g=c_{n}\int \frac{g(z)}{|x-z|^{n-1/2}}dz
\end{equation*}
(here and in the
following, integrals extend over all $\mathbb{R}^{n}$).
After integration by parts
we can split the quantity to estimate as follows:
\begin{equation*}
  \xx^{1/2+\epsilon}|D_{x}|^{-1/2}
  (\nabla_{x}\chi \cdot \nabla_{x}u)\simeq
  \int \beta(x,z)u(z)dz+
  \int \gamma(x,z)\nabla u(z)dz
\end{equation*}
where
\begin{equation*}
  \beta(x,z)=-\nabla_{z}\left(
    \frac{\xx^{1/2+\epsilon}\psi(x) [\nabla\chi(z)-\nabla\chi(x)]}
    {|x-z|^{n-1/2}}
  \right)
\end{equation*}
and
\begin{equation*}
  \gamma(x,z)=\frac{\xx^{1/2+\epsilon}\psi(x)\nabla\chi(x)}{|x-z|^{n-1/2}}.
\end{equation*}
In the following we extend the function $u$ as 0 outside $\Omega$
but keep the same notation for brevity. 
We have
\begin{equation*}
  \int \gamma(x,z)u(z)dz=
  \xx^{1/2+\epsilon}\psi(x)\nabla\chi(x)|D_{x}|^{-1/2}\nabla_{x} u
\end{equation*}
which implies, since $\psi$ has compact support,
\begin{equation*}
  \left\|\int \gamma(x,z)u(z)dz\right\|_{L^{2}_{x}}\lesssim
  \|\xx^{-1-\epsilon}|D_{x}|^{-1/2}\nabla_{x} u\|_{L^{2}_{x}}\lesssim
  \|\xx^{-1-\epsilon}|D_{x}|^{1/2} u\|_{L^{2}_{x}}
\end{equation*}
where in the last step we used \eqref{eq:equivw}. Finally, $\beta$
satisfies for all $N$
\begin{equation*}
  |\beta(x,z)|\lesssim\zz^{-N}|x-z|^{\frac12-n}
\end{equation*}
so that
\begin{equation*}
  \left\|\int \beta(x,z)u(z)dz\right\|_{L^{2}_{x}}\lesssim
  \||x|^{\frac12-n}*(\zz^{-N}u)\|_{L^{2}_{x}}\lesssim
  \|\zz^{-N}u\|_{L^{\frac{2n}{n+2}}}\lesssim
  \|\xx^{-1-\epsilon}u\|_{L^{2}_{x}}
\end{equation*}
by Hardy-Littlewood-Sobolev followed by H\"older's inequality
(for $N$ large enough). Summing up, and integrating also in the
remaining variables $t,y$, we arrive at
\begin{equation*}
  \|III\|_{L^{2}_{t}L^{2}_{y}L^{\frac{2n}{n-2}}_{x}}\lesssim
  \|\xx^{-1-\epsilon}u\|_{L^{2}_{t}L^{2}_{y}L^{2}_{x}(\Omega)}+
  \|\xx^{-1/2-\epsilon}|D_{x}|^{1/2}u\|_{L^{2}_{t}L^{2}_{y}L^{2}_{x}(\Omega)}
  \lesssim\|f\|_{L^{2}(\Omega)}
\end{equation*}
by \eqref{eq:smooforu}, \eqref{eq:smooforubis}.
In conclusion, putting together the estimates for $I$, $II$, $III$,
we obtain
\begin{equation}\label{eq:estw}
  \|w\|_{L^{2}_{t}L^{2}_{y}L^{\frac{2n}{n-2}}_{x}} \lesssim
  (1+\|\xx^{1+\epsilon} V\|_{L^{2}_{y}L^{n}_{x}})
    \|f\|_{L^{2}(\Omega)}.
\end{equation}

The remaining part $v=\chi(x)u$ can be estimated directly
via the Sobolev embedding
\begin{equation}\label{eq:sobol}
  \|g\|_{L^{\frac{2n}{n-2}}(A)}
  \lesssim\|\nabla g\|_{L^{2}(A)}
\end{equation}
which holds for any open set $A \subset \mathbb{R}^{n}$
(even unbounded)
and any $g\in H^{1}_{0}(\Omega)$, with a constant independent
of $A$. Then we have
\begin{equation}\label{eq:estv}
  \|\chi u\|_{L^{2}_{t}L^{2}_{y}L^{\frac{2n}{n-2}}_{x}}
  \lesssim
  \|u\nabla \chi\|
     _{L^{2}_{t}L^{2}_{x,y}}
  +\|\chi \nabla u\|
     _{L^{2}_{t}L^{2}_{x,y}}\lesssim
  \|f\|_{L^{2}(\Omega)}+\||D_{x}|^{1/2}f\|_{L^{2}(\Omega)}
\end{equation}
again by \eqref{eq:smooforu}.
Summing up \eqref{eq:estw} and \eqref{eq:estv}, we have proved
the following

\begin{theorem}\label{the:strichschro}
  Assume the domain 
  $\Omega \subseteq \mathbb{R}^{n}_{x}\times \mathbb{R}^{m}_{y}$,
  with $n\ge3$ and $m\ge1$, has a Lipschitz boundary,
  is repulsive w.r.to the $x$ variables and satisfies assumption
  \eqref{eq:assO}. Assume the potential $V(x,y)$ satisfies on
  $\Omega$ the inequalities
  \begin{equation}\label{eq:assVdelbis}
    V(x,y)\ge0,
    \qquad
    -\partial_{x}(|x|V(x,y))\ge0.
  \end{equation}
  and the operator $H=-\Delta_{x,y}+V(x,y)$ with Dirichlet boundary
  conditions is selfadjoint on $L^{2}(\Omega)$.
  Then the Schr\"odinger flow of $H$ satisfies the following 
  endpoint Strichartz estimate
  for all $f\in H^{1}_{0}(\Omega)$:
  \begin{equation}\label{eq:strichnonflat}
    \|e^{itH}f\|_{L^{2}_{t}L^{2}_{y}L^{\frac{2n}{n-2}}_{x}}
    \lesssim
    (1+\|\xx^{1+\epsilon} V\|_{L^{2}_{y}L^{n}_{x}})
    \Bigl(
      \|f\|_{L^{2}(\Omega)}+\||D_{x}|^{1/2}f\|_{L^{2}(\Omega)}
    \Bigr).
  \end{equation}
\end{theorem}

\begin{remark}\label{rem:nonhomg}
Proving \emph{nonhomogeneous} Strichartz estimates is more
difficult because of analytical technicalities. 
Recall that the solution to the nonhomogenous Schr\"odinger equation
\begin{equation}\label{eq:schrnonh}
  iu_{t}+H u=F(t,x,y),\qquad
  u(0,x,y)=f(x,y)
\end{equation}
on $\Omega$ can be represented as
\begin{equation*}
  u=e^{itH}f+i\int_{0}^{t}e^{i(t-s)H}F(s)ds;
\end{equation*}
we have already estimated the first term in Theorem \ref{the:strichschro},
and it remains to study the Duhamel operator
\begin{equation}\label{eq:duhamel}
  \int_{0}^{t}e^{i(t-s)H}F(s)ds.
\end{equation}
To this end, we introduce the norm
\begin{equation}\label{eq:H12}
  \|g\|_{H^{1/2,0}(\Omega)}=
  \|g\|_{L^{2}(\Omega)}+\||D_{x}|^{1/2}g\|_{L^{2}(\Omega)}
  \simeq
  \|\langle D_{x}\rangle^{1/2}g\|_{L^{2}(\Omega)} 
\end{equation}
and the corresponding Hilbert space $H^{1/2,0}(\Omega)$ defined
as the closure of $C^{\infty}_{c}(\Omega)$ in this norm.
Moreover we denote by $H^{-1/2,0}(\Omega)$ the dual of this space;
its norm can be characterized as
\begin{equation*}
  \|g\|_{H^{-1/2,0}(\Omega)}
  \simeq
  \|\langle D_{x}\rangle^{-1/2}g\|_{L^{2}(\Omega)}.
\end{equation*}
Then estimate \eqref{eq:strichnonflat} can be written
\begin{equation}\label{eq:strichhom}
   \|e^{itH}f\|_{L^{2}_{t}L^{2}_{y}L^{\frac{2n}{n-2}}_{x}}
    \lesssim
      \|f\|_{H^{1/2,0}(\Omega)},\qquad
      f\in H^{1}_{0}(\Omega).
\end{equation}
By interpolation with the conservation of energy
\begin{equation*}
   \|e^{itH}f\|_{L^{2}_{t}L^{2}(\Omega)}
    =\|f\|_{L^{2}(\Omega)}\le
      \|f\|_{H^{1/2,0}(\Omega)}
\end{equation*}
we obtain the full family of Strichartz estimates
\begin{equation}\label{eq:strfull}
   \|e^{itH}f\|_{L^{p}_{t}L^{2}_{y}L^{q}_{x}}
    \lesssim
      \|f\|_{H^{1/2,0}(\Omega)},\qquad
\end{equation}
for all \emph{admissible couples} $(p,q)$ of indices, i.e., such that
\begin{equation}\label{eq:admind}
  \frac n2=\frac2p+\frac nq,\quad
  2\le q\le \frac{2n}{n-2}.
\end{equation}
By duality, for any $F(t,x,y)\in L^{2}_{t}H^{1}_{0}(\Omega)$,
we have also
\begin{equation}\label{eq:dualstrich}
  \left\|
    \langle D_{x}\rangle^{-1/2}
    \int e^{-sH}F(s)ds
  \right\|_{L^{2}(\Omega)}\le C(V)
  \|F\|_{L^{p'}_{t}L^{2}_{y}L^{q'}_{x}}
\end{equation}
for $(p,q)$ admissible.
We also notice that estimates \eqref{eq:strfull} can be written
in the form
\begin{equation}\label{eq:strfullbis}
  \|e^{itH}\DDx^{-1}f\|_{L^{p}_{t}L^{2}_{y}L^{q}_{x}}
   \lesssim
     \|\DDx^{-1/2}f\|_{L^{2}(\Omega)},\qquad
  \frac n2=\frac2p+\frac nq,\quad
  2\le q\le \frac{2n}{n-2}.
\end{equation}
Now we can combine \eqref{eq:dualstrich} and 
\eqref{eq:strfullbis} to obtain
\begin{equation}\label{eq:almostnonh}
  \left\|
    \int e^{itH}\DDx^{-1}e^{-isH}F(s)ds
  \right\|_{L^{p}_{t}L^{2}_{y}L^{q}_{x}}\lesssim
  C(V)\|F\|_{L^{\widetilde{p}'}_{t}L^{2}_{y}L^{\widetilde{q}'}_{x}}.
\end{equation}
We can apply a standard trick and use the
Christ-Kiselev lemma as in \cite{KeelTao98-a}, which
permits to replace the integral over $\mathbb{R}$ with a truncated
integral over $[0,t]$, provided the indices satisfy the
additional condition $p>\widetilde{p}'$. This implies the
estimate
\begin{equation}\label{eq:almostnonh2}
  \left\|
    \int_{0}^{t} e^{itH}\DDx^{-1}e^{-isH}F(s)ds
  \right\|_{L^{p}_{t}L^{2}_{y}L^{q}_{x}}\lesssim
  C(V)\|F\|_{L^{\widetilde{p}'}_{t}L^{2}_{y}L^{\widetilde{q}'}_{x}}
\end{equation}
for all $(p,q)$ and $(\widetilde{p},\widetilde{q})$ admissible
such that $(p,\widetilde{p})\neq(2,2)$. To complete the proof we would need
an additional functional analytic assumption: 
\emph{the operator $\DDx$ commutes with the flow $e^{itH}$}; this
happens for instance when $V \equiv 0$. Then replacing $F$ with
$\DDx F$ in \eqref{eq:almostnonh2} we finally obtain
\begin{equation}\label{eq:almostnonh3}
  \left\|
    \int_{0}^{t} e^{i(t-s)H}F(s)ds
  \right\|_{L^{p}_{t}L^{2}_{y}L^{q}_{x}}\lesssim
  \|\DDx F\|_{L^{\widetilde{p}'}_{t}L^{2}_{y}L^{\widetilde{q}'}_{x}},
\end{equation}
i.e., the solution of \eqref{eq:schrnonh} satisfies
\begin{equation}\label{eq:strichlast}
  \|u\|_{L^{p}_{t}L^{2}_{y}L^{q}_{x}}\lesssim
  \|\DDx^{1/2}f\|_{L^{2}(\Omega)}+
  \|\DDx F\|_{L^{\widetilde{p}'}_{t}L^{2}_{y}L^{\widetilde{q}'}_{x}}
\end{equation}
for all admissible couples $(p,q)$ and $(\widetilde{p},\widetilde{q})$
with $(p,\widetilde{p})\neq(2,2)$.

\end{remark}

\begin{remark}\label{rem:NLS}
  In forthcoming works we shall apply the above Strichartz estimates
  to investigate the existence of small global solutions for
  nonlinear Schr\"odinger and wave equations on non flat waveguides.
\end{remark}


\end{document}